\documentclass{elsarticle}

\usepackage{amssymb}

\newtheorem{thm}{Theorem}
\newproof{pf}{Proof}

\lccode`\-=`\-
\usepackage[final,                     
color,notref,notcite]{showkeys}         

\usepackage{tikz}
\usepackage{pgfplots}
\usepackage{verbatim}

\journal{arXiv.org}

\newcommand{\grad}{\mathop{\rm grad}\nolimits}
\renewcommand{\div}{\mathop{\rm div}\nolimits}

\begin{document}

\begin{frontmatter}

\title{Substructuring domain decomposition scheme for unsteady problems}

\author{Petr N. Vabishchevich}
\ead{vab@ibrae.ac.ru}
\address{Nuclear Safety Institute, 
52, B. Tulskaya, 115191 Moscow, Russia}

\begin{abstract}

Domain decomposition methods are used for approximate solving boundary problems for partial differential equations on parallel 
computing systems. 
Specific  features of unsteady problems are taken into account in the most complete way in iteration-free schemes of domain decomposition. 
Regionally-additive schemes are based on different classes of splitting schemes. 
In this paper we highlight a class of domain decomposition schemes
which is based on the partition of the initial domain into subdomains with common  boundary nodes. 
Using the partition of unit we have constructed and studied unconditionally stable  schemes of domain decomposition
based on two-component splitting: the problem within subdomain and the problem at their boundaries. 
As an example there is considered the Cauchy problem for evolutionary equations 
of first and second order with non-negative self-adjoint operator in a finite Hilbert space. 
The theoretical consideration is supplemented with numerical solving a model problem 
for the two-dimensional parabolic equation. 
\end{abstract}

\begin{keyword}
unsteady problems \sep finite difference method 
\sep domain decomposition method 
\sep additive schemes \sep operator-splitting difference schemes

\MSC 65N06 \sep 65M06

\end{keyword}

\end{frontmatter}

\section{Introduction}
\label{sec:1}

Theory and practice of iterative solving stationary boundary value problems 
for partial differential equations is presented comprehensively in the books 
\cite{pre05281749,0931.65118,0857.65126,1069.65138}.
Different versions of the domain decomposition method with and without overlapping of subdomains are used. 
The approximate solution of unsteady problems can be derived via 
the standard implicit approximations in time and solving the corresponding 
grid problems at the new time level using one or another variants
of the domain decomposition method for stationary problems.  
Taking into account the transient character of unsteady problems (see, for example, the implementation 
on the basis of the Schwartz method \cite{Cai:1991:ASA,Cai:1994:MSM}), we can construct 
the optimal iterative methods of domain decomposition  
where the number of iterations is independent of the discretization steps in time and space. 

Specific  features of unsteady problems are taken into account in the most complete way in iteration-free schemes of domain decomposition. 
In some cases it is possible \cite{0825.65066,0766.65089} to employ only one iteration
of the Schwarz alternating method for the second order parabolic equation
without loss of accuracy of the approximate solution.  
Iteration-free domain decomposition schemes are associated with certain 
variants of the additive (splitting) schemes  --- regionally-additive schemes  \cite{1018.65103}.  

Domain decomposition schemes for solving  unsteady problems can be classified by 
the method of domain decomposition, choice of decomposition operators 
(exchange of boundary conditions) and used splitting scheme. 
For differential problems it is natural to select domain decomposition methods 
\begin{equation}\label{1}
  \overline{\Omega} = \bigcup_{\alpha =1}^{p} \overline{\Omega}_{\alpha},
  \quad \overline{\Omega}_{\alpha} = 
  \Omega_{\alpha} \cup \partial \Omega_{\alpha},
  \quad \alpha = 1, 2, ..., p  
\end{equation}
with overlapping of subdomains 
($\Omega_{\alpha \beta} \equiv  \Omega_{\alpha } \cap \Omega_{\beta} \neq \varnothing$)
and without overlapping ($\Omega_{\alpha \beta} = \varnothing$)
\cite{0931.65118,1069.65138}. 
Methods without overlapping of the subdomains are associated with an explicit formulation 
of boundary conditions at the interface boundaries. 
These methods are in common use for solving problems where in each particular
subdomain its own specific computational grid (triangulation) is introduced.  
To construct homogeneous computational algorithms,  
domain decomposition schemes with the overlapping of subdomains are employed.  
At the minimal overlapping where the width of overlapping is equal to the grid step 
($\Omega_{\alpha \beta} = \mathcal{O}(h)$), domain decomposition methods 
with the overlapping of subdomains can often be interpreted as methods 
without the overlapping of subdomains supplemented with appropriate boundary conditions of the exchange. 

Domain decomposition (\ref{1}) is associated with an appropriate additive 
representation of the problem  operator:
\begin{equation}\label{2}
  \mathcal{A} = \sum_{\alpha =1}^{p} \mathcal{A}_{\alpha} .
\end{equation}
In this case, the operator  $\mathcal{A}_{\alpha}$ is
associated with the solution of some problem in the subdomains
$\Omega_{\alpha}, \ \alpha = 1,2, ..., p$.
The most common approach to construct the operators of decomposition 
for solving boundary value problems for partial differential equations 
is based on using of the partition of unit for the computational domain. 
For decomposition  (\ref{1}) we can each separate subdomain 
$\Omega_{\alpha}$  associates with the function 
$\eta_{\alpha}(\mathbf{x}), \ \alpha = 1,2,...,p$ such that 
\begin{equation}\label{3}
  \eta_{\alpha}(\mathbf{x}) = \left \{
   \begin{array}{cc}
     > 0, &  \mathbf{x} \in \Omega_{\alpha},\\
     0, &  \mathbf{x} \notin  \Omega_{\alpha}, \\
   \end{array}
  \right .
  \quad \alpha = 1,2,...,p ,  
\end{equation}
and also  
\begin{equation}\label{4}
  \sum_{\alpha =1}^{p} \eta_{\alpha}(\mathbf{x}) = 1,
  \quad \mathbf{x} \in \Omega .
\end{equation}
Suppose, for example, that the operator $A$ is the diffusion operator: 
\begin{equation}\label{5}
  \mathcal{A} = - \div k(\mathbf{x}) \grad ,
  \quad \mathbf{x} \in \Omega .
\end{equation}
Then for the operators of decomposition we can define the following three basic forms: 
\begin{equation}\label{6}
  \mathcal{A}_{\alpha} =  \eta_{\alpha} \, \mathcal{A},
\end{equation}
\begin{equation}\label{7}
  \mathcal{A}_{\alpha} =  - \div k(\mathbf{x}) \eta_{\alpha} (\mathbf{x}) \grad ,
\end{equation}
\begin{equation}\label{8}
  \mathcal{A}_{\alpha} =  \mathcal{A} \, \eta_{\alpha} ,
  \quad \alpha = 1,2,...,p .
\end{equation}
This technique is used beginning with the work 
\cite{Laevsky} (decomposition  (\ref{7})), 
\cite{0719.65072} (decomposition (\ref{6})--(\ref{8})), 
the results of more recent works are summarized in the books \cite{1018.65103,0963.65091}. 
Various versions of the decomposition operators correspond to using  
different exchange boundary conditions and ensure the convergence of 
approximate solution in different spaces of grid functions.
Special attention should be given to issues of constructing decomposition operators 
for unsteady problems with non-selfadjoint operators \cite{0928.65102,vab_255,0888.65097}.

For unsteady problems with splitting (\ref{2}) different splitting schemes are used. 
In the theory of additive operator-difference schemes 
\cite{0963.65091,0971.65076,Marchuk:1990:SAD,0209.47103}
we need to distinguish the case of the simplest two-component splitting. 
In this case, we construct unconditionally stable factorized splitting schemes, 
such as the classical scheme of alternating directions, predictor-corrector scheme. 
Two-component regionally-additive schemes are constructed and studied in 
\cite{Laevsky,0719.65072,0838.65086} as well as in the above papers 
\cite{0928.65102,vab_255,0888.65097}
for convection-diffusion problems. 
 
In application of domain decomposition methods
the splitting of problem operator into the sum of three or more
non-commutative operators ($p > 2$ in (\ref{2})) is of great interest.
Classic schemes \cite{0971.65076,Marchuk:1990:SAD,0209.47103}
of multi-component splitting are based on the concept of summarized approximation. 
Additively-averaged schemes of summarized approximation \cite{0963.65091,0297.35037}
are more explicitly oriented to parallel computations. 
Regionally-additive schemes of component-wise splitting are investigated in  \cite{0986.65510}. 
A variant of two-component splitting with the Crank-Nicolson scheme for the 
individual subproblems with the minimal overlapping and decomposition (\ref{7})  
is considered in the article  \cite{dryja2007}.

Nowadays, the schemes of full approximation are in common use for the general multi-component splitting. 
In this regard, we note regularized additive schemes \cite{samarskii1998regularized}
where the condition of stability is achieved due to perturbations of operators of the difference scheme. 
In the vector additive schemes \cite{0712.65089,vabishchevich1996vector}
instead of one equation we solve a system of similar equations. 
Such schemes are also constructed for the evolutionary equations of second order \cite{samarskii1992regularized,abrashin1998numerical}.
Vector regionally-additive schemes are investigated in  \cite{0863.65056,0965.65119}. 
In the work \cite{1156.65084} there are proposed more general regularized schemes of domain decomposition  
with different structures for both the splitting operators and operators of the grid problem at a new time-level. 

Among other domain decomposition methods for solving boundary 
value problems for parabolic equations it is necessary to highlight  explicit-implicit 
methods considered in many papers (see, for example, \cite{1013.65106,1120.65098,1076.65079,1099.65084,1094.65097,1168.65387}).
The domain decomposition in this case is performed without the overlapping of subdomains
and the transition to a new time-level is organized as follows. 
First, the approximate solution at the common boundaries of subdomains is predicted using the explicit scheme.
Next, these boundary conditions are used to derive the approximate solution within individual subdomains.
And finally, a correction of interface boundary conditions is carried out using implicit schemes. 
It will be shown below, that such schemes of domain decomposition 
are completely fit in the above general scheme of decomposition methods 
at a special domain decomposition with the choice of operators according to 
decomposition (\ref{6}).

In this paper we construct domain decomposition schemes
for parabolic and hyperbolic equations with self-adjoint elliptic operators of second order. 
Unconditionally stable factorized regionally-additive schemes are  constructed
using decomposition  (\ref{6}) and the two-component and general multi-component splitting. 
Domain decomposition schemes with a self-adjoint operator for the grid problem at
the new time-level are derived, that allows to construct on their basis iterative methods of domain 
decomposition--conjugate gradients for solving stationary problems. The paper is organized as follows. 
In Section 2 we formulate the model Cauchy problem for parabolic 
and hyperbolic equations in a rectangle. 
Next, Section 3 provides the stability conditions for the standard two- and tree-level implicit 
schemes with weights for model problems in a rectangle. 
The domain decomposition and construction of the operators are discussed in Section 4.  
The possibilities of the standard factorized schemes for domain decomposition are considered in Section 5. 
The stability condition, appropriate a priori estimates for 
the approximate solution and convergence rate estimate are derived for factorized regionally-additive schemes. 
In Section 6 there are constructed multi-component splitting schemes. 
Hyperbolic equations of second order are considered in Section 7. 
Theoretical results are illustrated by the numerical results presented in Section 8. 

\section{Model boundary problems}
\label{sec:2}

Let us consider a model boundary value problem for the parabolic equation of second order.  
In a bounded domain $\Omega$ the unknown function  $u(\mathbf{x},t)$ satisfies the following equation
\begin{equation}\label{9}
   \frac{\partial u}{\partial t} 
   - \sum_{\alpha =1}^{m}
   \frac{\partial }{\partial x_\alpha} 
   \left ( k({\bf x})  \frac{\partial u}{\partial x_\alpha} \right ) = f({\bf x},t),
   \quad {\bf x}\in \Omega,
   \quad 0 < t \leq T,
\end{equation}
where  $k(\mathbf{x}) \geq \kappa > 0, \  {\bf x}\in \Omega$.
Equation  (\ref{9}) is supplemented with the homogeneous Dirichlet boundary conditions 
\begin{equation}\label{10}
   u({\bf x},t) = 0,
   \quad {\bf x}\in \partial \Omega,
   \quad 0 < t < T.
\end{equation}
In addition, the initial condition is prescribed 
\begin{equation}\label{11}
   u({\bf x},0) = u^0({\bf x}),
   \quad {\bf x}\in \Omega.
\end{equation}

Unsteady diffusion problem (\ref{9})--(\ref{11}) is considered 
on the set of functions  $u({\bf x},t)$ 
satisfying boundary conditions (\ref{10}). 
Then instead of  (\ref{9}), (\ref{10}) we use the operator-differential equation 
\begin{equation}\label{12}
  \frac {du}{dt} + {\cal A} u = f(t),
  \quad 0 < t < T.
\end{equation}
The  Cauchy problem is considered for evolutionary equation  (\ref{12}):
\begin{equation}\label{13}
  u(0) = u^0 .
\end{equation}
For the diffusion operator we assume  
\[
  {\cal A} u = 
   - \sum_{\alpha =1}^{m}
   \frac{\partial }{\partial x_\alpha} 
   \left ( k({\bf x})  \frac{\partial u}{\partial x_\alpha} \right ) .
\]

On the set of functions (\ref{10}) let us define the Hilbert space ${\cal H} =  {\cal L}_2 (\Omega)$ 
with the scalar product and norm 
\[
  (u, v) = \int_{\Omega} u(\mathbf{x}) v(\mathbf{x}) d \mathbf{x},
  \quad \| u\| = (u,u)^{1/2} .
\]
In ${\cal H}$ the operator of the diffusive transport ${\cal A}$ is self-adjoint and positive definite: 
\begin{equation}\label{14}
  {\cal A} = {\cal A}^* \ge \kappa \delta {\cal E},
  \quad \delta  = \delta(\Omega) > 0,
\end{equation}
where ${\cal E}$ is  the identity operator in  ${\cal H}$.

We present now the simplest a priori estimate for the solution of problem (\ref{12})--(\ref{14})
which will be for us the check point for the considering grid problems. 
The self-adjoint positive definite operator $\mathcal{D}$ can be
associated with the Hilbert space $\mathcal{H}_{\mathcal{D}}$ having the inner product and norm
\[
  (u,v)_{\mathcal{D}} = (\mathcal{D}u,v),
  \quad \| u\|_{\mathcal{D}} = (u,u)^{1/2}_{\mathcal{D}}
\]
respectively.
In $\mathcal{H}$ multiply scalarly equation (\ref{12}) by $\mathcal{A} u$. 
In view of  (\ref{14}) we obtain inequality
\begin{equation}\label{15}
  \frac{1}{2} \frac{d}{d t} \|u\|_{\mathcal{A}}^2 +
  \| \mathcal{A} u \|^2 = (f,\mathcal{A} u) .
\end{equation}
Taking into account 
\[
  (f, \mathcal{A} u) \leq \| \mathcal{A} u \|^2
  + \frac{1}{4} \| f\|^2 ,
\]
from (\ref{15}) we have 
\[
  \frac{d}{d t} \|u\|_{\mathcal{A}}^2 \leq \frac{1}{4} \| f\|^2 .
\]
In view of the Gronwall lemma we obtain the desired estimate 
\begin{equation}\label{16}
  \|u\|_{\mathcal{A}}^2 \leq \|u^0\|_{\mathcal{A}}^2 + 
  \int_{0}^{t} \|f(\theta)\|^2 d \theta ,
\end{equation}
which expresses the stability of the solution of problem  (\ref{12})--(\ref{14})  
with respect to the initial data and right-hand side. 

In addition to the parabolic equation (\ref{9}), we consider the hyperbolic equation 
\begin{equation}\label{17}
   \frac{\partial^2 u}{\partial t^2} 
   - \sum_{\alpha =1}^{m}
   \frac{\partial }{\partial x_\alpha} 
   \left ( k({\bf x})  \frac{\partial u}{\partial x_\alpha} \right ) = f({\bf x},t),
   \quad {\bf x}\in \Omega,
   \quad 0 < t \leq T
\end{equation}
with boundary conditions  (\ref{10}).
Equation (\ref{17})  is supplemented with two initial conditions 
\begin{equation}\label{18}
   u({\bf x},0) = u^0({\bf x}),
   \quad \frac{\partial u}{\partial t} ({\bf x},0) = v^0({\bf x}),
   \quad {\bf x}\in \Omega.
\end{equation}
Problem  (\ref{10}), (\ref{17}), (\ref{18})  is associated with the following 
Cauchy problem for the  evolutionary equation of second order: 
\begin{equation}\label{19}
  \frac {d^2 u}{dt^2} + {\cal A} u = f(t),
  \quad 0 < t < T,
\end{equation}
\begin{equation}\label{20}
  u(0) = u^0 ,
  \quad \frac {d u}{dt}(0) = v^0 .
\end{equation}

Multiply scalarly equation (\ref{19})  by $\mathcal{A} du/dt$ and obtain 
\[
  \frac{1}{2} \frac{d}{dt} 
  \left (\left \|\frac{du}{dt} \right \|_{\mathcal{A}}^{2} + 
  \|\mathcal{A} u\|^2 \right ) =
  \left (f, \mathcal{A} \frac{d u}{dt} \right ).
\]
For the right-hand side we use the estimate
\[
  \left (f, \mathcal{A} \frac{d u}{dt} \right ) \leq
  \frac{1}{2} \left \|\frac{du}{dt} \right \|_{\mathcal{A}}^2 + 
  \frac{1}{2} \|f\|_{\mathcal{A}}^2.
\]
The result is 
\[
  \frac{d}{dt} \|u\|_*^2  \leq \|u\|_*^2 + \|f\|_{\mathcal{A}}^2,
\]
where 
\[
  \|u\|_*^2 = \left \|\frac{du}{dt} \right \|_{\mathcal{A}}^{2} + 
  \|\mathcal{A} u\|^2.
\]
The desired a priori estimate 
\begin{equation}\label{21}
  \|u(t)\|_*^2 \leq \exp (t) 
  \left (\|\mathcal{A} u^0\|^2 + \|v^0\|_{\mathcal{A}}^2 +
  \int\limits_{0}^{t}{\exp (-\theta) \|f(\theta)\|_{\mathcal{A}}^2} d \theta \right )
\end{equation}
expresses the stability with respect to the initial data and right-hand side of the Cauchy problem 
for operator-differential equation (\ref{19}).

\section{Standard difference approximations}
\label{sec:3}

We will conduct a detailed study of approximations in space and time using as 
an example the boundary problems in a rectangle 
\[
  \Omega = \{ \ \mathbf{x} \ | \ \mathbf{x} = (x_1, x_2), 
  \ 0 < x_{\alpha} < l_{\alpha}, \ \alpha =1,2 \}.
\]
The approximate solution is given at the nodes of a uniform rectangular grid $\Omega$:
\[
   \bar{\omega} = \{ \mathbf{x} \ | \ \mathbf{x} = (x_1, x_2),
   \quad x_\alpha = i_\alpha h_\alpha,
   \quad i_\alpha = 0,1,...,N_\alpha,
   \quad N_\alpha h_\alpha = l_\alpha\} 
\]
and let $\omega$ be the set of internal nodes 
($\bar{\omega} = \omega \cup \partial \omega$). 
For the grid functions $y(\mathbf{x}) = 0, \ \mathbf{x} \in \partial \omega$
we define the Hilbert space $H = L_2({\omega})$ 
with the scalar product and norm 
\[
  (y,w) = \sum_{{\bf x} \in \omega}
  y({\bf x}) w({\bf x}) h_1 h_2,
  \quad \|y\| = (y,y)^{1/2} .
\]

Assuming that the coefficient  $k(\mathbf{x})$ in $\Omega$  is sufficiently smooth, 
we take the grid operator of the diffusion as 
\[
  A y =
  - \frac{1}{h_1^2} k(x_1+0.5h_1,x_2)
  (y(x_1+h_1,x_2) - y(x_1,x_2))
\]
\[
  + \frac{1}{h_1^2} k(x_1-0.5h_1,x_2)
  (y(x_1,x_2) - y(x_1-h_1,x_2)) 
\]
\[
  - \frac{1}{h_2^2} k(x_1,x_2+0.5h_2)
  (y(x_1,x_2+h_2) - y(x_1,x_2))
\]
\begin{equation}\label{22}
  + \frac{1}{h_2^2} k(x_1,x_2-0.5h_2)
  (y(x_1,x_2) - y(x_1,x_2-h_2)) .
\end{equation}
In $H$  the operator  $A$ is self-adjoint and positive definite:
\begin{equation}\label{23}
  A = A^* \geq \kappa (\delta_1+\delta_2) E,
  \quad \delta_{\alpha} = 
  \frac{4}{h^2_{\alpha}} \sin^2 \frac{\pi h_{\alpha}}{2 l_{\alpha}} ,
  \quad \alpha = 1,2.
\end{equation}

After the approximation in space we go from (\ref{9}), (\ref{10}) 
to the differential-difference equation 
\begin{equation}\label{24}
  \frac{d y}{d t} + A y = f(\mathbf{x},t),
  \quad \mathbf{x} \in \omega,
   \quad 0 < t < T .
\end{equation}
Taking into account  (\ref{11}), let us supplement equation (\ref{24}) with the initial condition 
\begin{equation}\label{25}
   y({\bf x},0) = u^0({\bf x}),
   \quad {\bf x}\in \omega .
\end{equation}
For the solution of  the differential-difference Cauchy problem (\ref{24}), (\ref{25}) the following a priori estimate holds
(see (\ref{16}))
\begin{equation}\label{26}
  \|y\|_{A}^2 \leq \|u^0\|_{A}^2 + 
  \int_{0}^{t} \|f(\theta)\|^2 d \theta .
\end{equation}

Similarly, the approximation in space leads us from (\ref{10}), (\ref{17}), (\ref{18}) to the problem 
\begin{equation}\label{27}
  \frac{d^2 y}{d t^2} + A y = f(\mathbf{x},t),
  \quad \mathbf{x} \in \omega,
   \quad 0 < t < T ,
\end{equation}
\begin{equation}\label{28}
   y({\bf x},0) = u^0({\bf x}),
   \quad \frac{d y}{d t}({\bf x},0) = v^0({\bf x}),
   \quad {\bf x}\in \omega .
\end{equation}
The grid analog of (\ref{21}) is the estimate 
\begin{equation}\label{29}
  \|y(t)\|_*^2 \leq \exp (t) \left (\|A u^0\|^2 + \|v^0\|_{A}^2 +
  \int\limits_{0}^{t}{\exp (-\theta) \|f(\theta)\|_{A}^2} d \theta \right ),
\end{equation}
where 
\[
  \|y\|_*^2 = \left \|\frac{dy}{dt} \right \|_{A}^{2} + \|A y\|^2.
\]

The emphasis now is on the approximation in time. 
In the construction of domain decomposition schemes for problem  (\ref{24}), (\ref{25}), 
the starting point for us is the usual two-level schemes. 
Let $\tau$  be of a uniform time-step and let 
 $y^n = y(t^n), \ t^n = n \tau$, $n = 0,1, ..., N, \ N\tau = T$.
Equation  (\ref{24}) is approximated by a two-level scheme with weights 
\begin{equation}\label{30}
  \frac{y^{n+1} - y^{n}}{\tau }
  + A(\sigma y^{n+1} + (1-\sigma) y^{n}) = \varphi^n,
  \quad n = 0,1, ..., N-1,
\end{equation}
where, for example, $\varphi^n = f(\sigma t^{n+1} + (1-\sigma) t^{n})$.
It is supplemented by the initial condition 
\begin{equation}\label{31}
  y^0 = u^0 .
\end{equation}
Difference scheme  (\ref{30}), (\ref{31}) has the approximation error 
$\mathcal{O} (\tau^2 + (\sigma - 1/2) \tau + h^2)$, where 
$h^2 = (h_1^2 + h_2^2)/2$. 

\begin{thm} 
\label{t-1} 
Difference scheme (\ref{30}), (\ref{31}) is unconditionally stable for 
$\sigma \geq 1/2$,
and for the numerical solution the estimate 
\begin{equation}\label{32}
  \|y^{n+1}\|^2_{D} \leq \|y^{n}\|^2_{D} + 
  \frac{\tau}{2} \|\varphi^n\|^2,
  \quad n = 0,1, ..., N-1 ,
\end{equation}
holds, where 
\[
  D = A + \left (\sigma - \frac{1}{2}\right )  \tau A^2 .
\]
\end{thm} 

\begin{pf} 
Let write difference scheme (\ref{30})  as 
\[
  \left (E + \left (\sigma - \frac{1}{2} \right ) \tau A \right )
  \frac{y^{n+1} - y^{n}}{\tau }
  + A \frac{y^{n+1} + y^{n}}{2 } = \varphi^n,
\]
and  multiply scalarly it by $\tau A (y^{n+1} + y^{n})$. 
Using the fact that $\sigma \geq 1/2$ the operator $D \geq A$, we have 
\[
  \|y^{n+1}\|^2_{D} - \|y^{n}\|^2_{D} +
  \frac{\tau}{2} \| A(y^{n+1} + y^{n}) \|^2 =
  \tau (\varphi^n, A(y^{n+1} + y^{n})).
\]
Taking into account 
\[
  (\varphi^n, A(y^{n+1} + y^{n})) \leq 
  \frac{1}{2} \| A(y^{n+1} + y^{n}) \|^2 +
  \frac{1}{2} \|\varphi^n\|^2,
\]
we obtain the required estimate (\ref{32}).
\end{pf} 

A priori estimate  (\ref{32}) for the solution of problem (\ref{30}), (\ref{31}) is a 
grid analog of the a priori estimate (\ref{26}) for the solution of differential-difference problem 
(\ref{24}), (\ref{25}) ($D = A + \mathcal{O}(\tau)$).

To solve numerically problem (\ref{27}), (\ref{28}), it is natural to use three-level schemes of second order accuracy in time. 
Let 
\[
  \frac{y^{n+1} - 2y^{n} + y^{n-1}}{\tau^2 }
  + A(\sigma y^{n+1} + (1-2 \sigma) y^{n} + \sigma y^{n-1}) = \varphi^n,
\]
\begin{equation}\label{33}
  n = 1,2, ..., N-1,
\end{equation}
where, for example,  $\varphi^n = f(t^{n})$.
In view of (\ref{28}) we can for the solution of equation (\ref{27}) approximate the initial condition as follows: 
\begin{equation}\label{34}
  y^0 = u^0 ,
  \quad \frac{y^1 - y^0}{\tau} = v^0 +
  \frac{\tau}{2} (\varphi^0 - A u^0) .
\end{equation}
The error of difference scheme (\ref{33}), (\ref{34}) is 
$\mathcal{O} (\tau^2 + h^2)$.
 
\begin{thm} 
\label{t-2} 
Difference scheme  (\ref{33}), (\ref{34}) is unconditionally stable for $\sigma \geq 1/4$,
and for the numerical solution the estimate 
\begin{equation}\label{35}
  S^{n+1} \leq \exp(\tau) S^{n} 
  + \frac{\tau^2}{2} \frac{\exp(\tau)}{\exp(0.5\tau)-1}  \|\varphi^n \|^2_A,
  \quad n = 0,1, ..., N-1 ,
\end{equation}
holds, where 
\[
  S^n =  \left \|\frac{y^n - y^{n-1}}{\tau } \right \|^2_D + 
  \left \|A \frac{y^n + y^{n-1}}{2} \right \|^2 ,
\]
\[
  D = A + \left (\sigma - \frac{1}{4} \right ) \tau^2 A^2 .
\]
\end{thm} 

\begin{pf} 
We introduce the notation 
\[
  \zeta^n = \frac{y^n + y^{n-1}}{2} ,
  \quad \eta^n = \frac{y^n - y^{n-1}}{\tau } .
\]
Taking into account the identity 
\[
  y^n = \frac{1}{4} ( y^{n+1} + 2 y^n + y^{n-1})
  - \frac{1}{4} ( y^{n+1} - 2 y^n + y^{n-1}),
\]
\[
  \sigma y^{n+1} + (1- 2\sigma) y^{n} + \sigma y^{n-1} = 
  y^n  + \sigma ( y^{n+1} - 2 y^n + y^{n-1})
\]
we rewrite  (\ref{33}) as 
\begin{equation}\label{36}
  \left (E + \left (\sigma - \frac{1}{4} \right ) \tau^2 A \right ) \frac{\eta^{n+1} - \eta^{n}}{\tau} + 
  A \frac{\zeta^{n+1} + \zeta^{n}}{2} = \varphi^n .
\end{equation}
Multiply scalarly (\ref{36}) in  $H$ by 
\[
  2 A(\zeta^{n+1} - \zeta^{n}) = \tau A (\eta^{n+1} + \eta^{n}) .
\]
With this notation for $\sigma \geq 1/4$ we obtain 
\begin{equation}\label{37}
  S^{n+1} - S^{n} =
  \tau (\varphi^n, A(\eta^{n+1} + \eta^{n})) .
\end{equation}
Using the estimates for the right-hand side 
\[
  \tau A (\varphi^n,(\eta^{n+1} + \eta^{n})) \leq 
  \frac{\tau  }{2 \varepsilon}   \|\eta^{n+1} + \eta^{n}\|^2_A +
  \frac{\tau}{2}  \varepsilon  \|\varphi^n \|^2_A ,
\]
\[
  \|\eta^{n+1} + \eta^{n}\|^2_A \leq 
  2 (\|\eta^{n+1}\|^2_A  + \|\eta^{n}\|^2_A ),
\]
with $\varepsilon > 0$, from (\ref{37}) we obtain 
\begin{equation}\label{38}
  \left (1 - \frac{\tau }{\varepsilon } \right )  S^{n+1} \leq 
  \left (1 + \frac{\tau }{\varepsilon } \right )  S^{n} +
  \frac{\tau}{2}  \varepsilon \|\varphi^n \|^2_A .
\end{equation}
We choose $\varepsilon$ so that 
\[
  1 - \frac{\tau }{\varepsilon } = \exp(-0.5 \tau) , 
\] 
and therefore 
\[
  1 + \frac{\tau }{\varepsilon } = \exp(0.5 \tau) .
\]
With this in mind from (\ref{38}) we obtain the level-wise stability estimate (\ref{35}). 
\end{pf} 

Estimate (\ref{35}) can be treated as the grid analog of the a priori estimate (\ref{29}). 
For difference schemes (\ref{30}), (\ref{31}) and  (\ref{33}), (\ref{34})
we can obtain many other a priori estimates of stability with respect to the initial data and right-hand side 
\cite{1018.65103,0971.65076}. 
We have restricted to only those estimates that we can associate with
the corresponding estimates for the domain decomposition schemes considered below.

\section{Substructuring domain decomposition}
\label{sec:4}

\begin{figure}[ht] 
  \begin{center}
	\begin{tikzpicture}[domain=0:4]
		\draw[step=0.5,very thin,color=gray] (0,0) grid (8,8);
		\draw[step=4,color=black,line width=1pt] (0,0) grid (8,8);
		\draw[|<->|] (8.25,2) -- (8.25,2.5);
		\foreach \x in {1,...,15}
  			\foreach \y in {1,...,15}
	  			\draw [{fill=white}] (0.5*\x,0.5*\y) circle (0.1);
		\foreach \x in {1,...,15}
	  		\fill [black] (0.5*\x,4) circle (0.1);
		\foreach \y in {1,...,15}
	  		\fill [black] (4,0.5*\y) circle (0.1);
		\draw(8.5,2.3) node {$h$};
		\draw[|<->|] (8.25,4) -- (8.25,8);
		\draw(8.5,6) node {$\widehat{h}$};
	\end{tikzpicture}
    \caption{Grid decomposition} 
	\label{f-1}
  \end{center}
\end{figure}
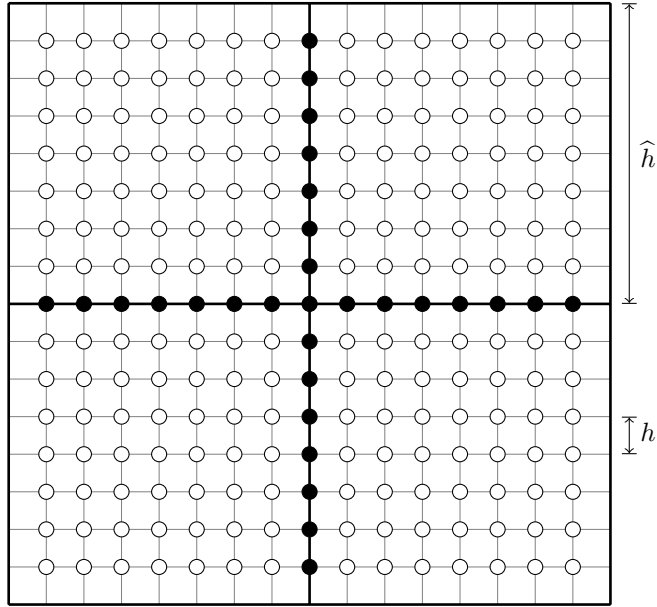

\begin{figure}[ht] 
  \begin{center}
	\begin{tikzpicture}[domain=0:4]
		\fill [black!15] (0,0) rectangle +(8,8);
		\draw [black] (0,0) rectangle +(8,8);
		\fill [black!45] (3.75,0) rectangle +(0.5,8);
		\fill [black!45] (0,3.75) rectangle +(8,0.5);
 		\draw[step=0.5,very thin,color=black!0] (0,0) grid (8,8);	
		\draw(4,4) node {$\Omega_2$};
		\draw[|<->|] (8.25,3.75) -- (8.25,4.25);
		\draw(8.5,4.) node {$h$};
		\draw[|<->|] (3.75,8.25) -- (4.25,8.25);
		\draw(4,8.5) node {$h$};
		\draw(6,6) node {$\Omega_1$};
	\end{tikzpicture}
    \caption{Domain decomposition} 
	\label{f-2}
  \end{center}
\end{figure}
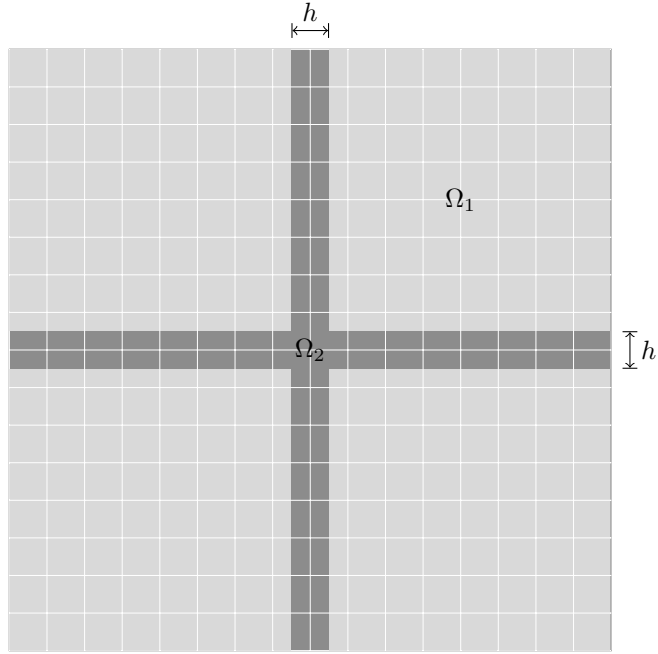

Let us consider a special class of domain decomposition methods. 
At the discrete level we define a set of interface nodes inside the domain
and then solve the subproblems separately  inside the subdomains. 
At the continuous level, this decomposition is associated with subdomains the width of which is equal to 
the corresponding  discretization step in space. 
We illustrate our consideration on the model grid problems in a rectangular. 

The computational grid  $\omega$ is partitioned into rectangular 
subdomains of a coarse grid with the step $\widehat{h}$. 
The boundaries of the subdomains (direct lines) consist of the nodes of the fine computational grid.  
Denote this set of interior boundary nodes as  $\widehat{\omega}$.
A fragment of the grid is shown in Fig.~\ref{f-1}.
Such a decomposition of the fine computational grid corresponds to 
the domain decomposition depicted in Fig.~\ref{f-2}:
$\overline{\Omega} = \overline{\Omega}_1 \cup \overline{\Omega}_2$, 
$\Omega_{12} = \varnothing$.
Subdomain $\Omega_2$ is a lattice, width of the individual edges of the lattice is  $h$.
Domain  $\Omega_1$  consists of disconnected individual subdomains. 
 
The partition of unit for  (\ref{3}), (\ref{4})  we associate with the corresponding additive 
representation of the identity operator  $E$  in the space of grid functions $H$, 
defined on the set of internal nodes of  $\omega$.
Let
\begin{equation}\label{39}
  \sum_{\alpha =1}^{p} \chi_{\alpha} = E,
  \quad \chi_{\alpha} \geq 0,
  \quad \alpha = 1,2,...,p .
\end{equation}
Similarly (\ref{6}), the operators of decomposition can be given in the form 
\begin{equation}\label{40}
  A_{\alpha} = \chi_{\alpha} A,
  \quad \alpha = 1,2,...,p .
\end{equation}
In view of  (\ref{39}), in this splitting we have for the problem operator
the following additive representation 
\begin{equation}\label{41}
  A = \sum_{\alpha =1}^{p}A_{\alpha} .
\end{equation}

Splitting (\ref{41}) allows us to go from equation (\ref{24}) to the equation 
\begin{equation}\label{42}
  \frac{d y}{d t} + \sum_{\alpha =1}^{p}A_{\alpha} y = f(\mathbf{x},t),
  \quad \mathbf{x} \in \omega,
   \quad 0 < t < T .
\end{equation}
Direct construction of various splitting schemes for problem  (\ref{25}), (\ref{42}) 
is complicated by the fact that individual operator terms $A_{\alpha}, \ \alpha = 1,2,...,p$ 
do not inherit the basic properties of the operator  $A$ --- the self-adjointness and non-negativity. 
However, using decomposition operators (\ref{40}), 
equation (\ref{42}) can be easy transformed in the symmetric form. 
Multiplying equation (\ref{42}) by the self-adjoint operator  $A$, we obtain the equation 
\begin{equation}\label{43}
  \tilde{B} \frac{d y}{d t} + \sum_{\alpha =1}^{p} \tilde{A}_{\alpha} y = A f(\mathbf{x},t),
  \quad \mathbf{x} \in \omega,
   \quad 0 < t < T ,
\end{equation}
where the operators 
\[
  \tilde{B} = A,
  \quad \tilde{A}_{\alpha} = A \chi_{\alpha} A,
  \quad \alpha = 1,2,...,p 
\]
are self-adjoint and non-negative. 
Moreover, we can introduce new variables $v = A^{1/2} y$ 
and instead of (\ref{43})  we can consider the equation 
\begin{equation}\label{44}
  \frac{d v}{d t} + \sum_{\alpha =1}^{p} \tilde{A}_{\alpha} v = A^{1/2} f(\mathbf{x},t),
  \quad \mathbf{x} \in \omega,
   \quad 0 < t < T ,
\end{equation}
with self-adjoint and non-negative operators 
\[
  \tilde{A}_{\alpha} = A^{1/2} \chi_{\alpha} A^{1/2},
  \quad \alpha = 1,2,...,p .
\]
Standard estimates for the solution of equation  (\ref{44}) 
in the norm of $H$ (for $\| v\|$) correspond to using
estimates in  $H_A$ (for $\|y\|_A$). 
This explains our unusual in some sense choice of the the priori estimate (\ref{26}) 
for problem (\ref{24}), (\ref{25}) and estimate (\ref{29}) for problem  (\ref{27}), (\ref{28}).

The particular specification of the decomposition operators of type  (\ref{39}), (\ref{40})
is provided via the selection of terms $\chi_{\alpha}, \ \alpha = 1,2,...,p$.
Some advanced features are discussed below, but we start from the simplest version. 
If we use substructuring domain decomposition (see Fig.~\ref{f-1}), it is natural to put 
\begin{equation}\label{45}
  \chi_{2}(\mathbf{x}) = \left \{
  \begin{array}{ll}
    1, &  \mathbf{x} \in \widehat{\omega}, \\
    0, &  \mathbf{x} \notin  \widehat{\omega}, \\
  \end{array}
  \right .
  \quad \chi_{1}(\mathbf{x}) = 1 - \chi_{2}(\mathbf{x}),
  \quad \mathbf{x} \in \omega .
\end{equation}
The operator $A_2$ is associated with interface nodes $\widehat{\omega}$,
whereas $A_1$ --- with the internal nodes of subdomains. 

\section{Factorized schemes of domain decomposition}
\label{sec:5}

After selecting the operators in decomposition (\ref{41}) 
the construction of domain decomposition schemes is carried out using one or another additive schemes. 
For (\ref{40}), (\ref{45}) we can consider the simplest two-component ($p=2$) splitting schemes. 
In this situation, we can try to use the operator analogues of the classical schemes of alternating directions\cite{pearac,dourac}.

We begin with the scheme of stabilizing correction \cite{dourac}, 
where the transition to a new time level in problem (\ref{25}), (\ref{42})
with $p=2$ is performed  as follows: 
\begin{equation}\label{46}
  \frac{y^{n+1/2} - y^{n}}{\tau }
  + A_{1} y^{n+1/2} + A_{2} y^{n} = \varphi^n,
\end{equation}
\begin{equation}\label{47}
  \frac{y^{n+1} - y^{n}}{\tau }
  + A_{1} y^{n+1/2} + A_{2} y^{n+1} = \varphi^n,
\end{equation}
where, for example,  $\varphi^n = f(t^{n+1})$, $n = 0,1, ..., N-1$.
Equations (\ref{46}), (\ref{47}) are complemented by the initial condition  (\ref{31}).

If decomposition  (\ref{40}), (\ref{45}) is used, we have 
\begin{equation}\label{48}
  \frac{y^{n+1/2} - y^{n}}{\tau }
  + \chi_{1} A y^{n+1/2} + \chi_{2} A  y^{n} = \varphi^n,
\end{equation}
\begin{equation}\label{49}
  \frac{y^{n+1} - y^{n}}{\tau }
  + \chi_{1} A  y^{n+1/2} + \chi_{2} A  y^{n+1} = \varphi^n .
\end{equation}
The implementation of this scheme can be different. 

Taking into account that 
\[
  \varphi^n = \chi_{1}\varphi^n + \chi_{2} \varphi^n,
\]
let us introduce the auxiliary function  $\widetilde{y}^{n+1/2}$ 
and divide equation (\ref{48}) into two ones: 
\begin{equation}\label{50}
  \frac{\widetilde{y}^{n+1/2} - y^{n}}{\tau }
  + \chi_{2} A  y^{n} = \chi_{2} \varphi^n,
\end{equation}
\begin{equation}\label{51}
  \frac{y^{n+1/2} - \widetilde{y}^{n+1/2}}{\tau }
  + \chi_{1} A  y^{n+1/2} = \chi_{1} \varphi^n .
\end{equation}
The function $\widetilde{y}^{n+1/2}$ is determined via the explicit scheme  (\ref{50}).
Moreover, taking into account (\ref{45}),
the calculations are performed only on the set of interface nodes. 

\textit{Stage 1}. Evaluation of the conditions at the boundaries of the subdomains via the explicit scheme:
\[
  \frac{\widetilde{y}^{n+1/2} - y^{n}}{\tau }
  + A  y^{n} = \varphi^n,
  \quad \mathbf{x} \in \widehat{\omega},
\]
\[
  \widetilde{y}^{n+1/2} = y^{n},
  \quad \mathbf{x} \notin  \widehat{\omega} .
\]
After such a predictor of boundary conditions we solve problems in subdomains  (\ref{51}). 

\textit{Stage 2}. Evaluation of the solution in subdomains using the implicit scheme:
\[
  \frac{y^{n+1/2} - \widetilde{y}^{n+1/2}}{\tau }
  + A  y^{n+1/2} = \varphi^n ,
  \quad \mathbf{x} \notin  \widehat{\omega} .
\]
\[
  y^{n+1/2} = \widetilde{y}^{n+1/2} ,
  \quad \mathbf{x} \in \widehat{\omega} .
\]

The last step is to correct conditions at the boundaries, which provides, 
in particular, the stability of the approximate solution. 
For the subdomains it is convenient to replace equation (\ref{49}) by the difference  of (\ref{49}), (\ref{48}):
\[
  \frac{y^{n+1} - y^{n+1/2}}{\tau }
  + \chi_{2} A  (y^{n+1} - y^{n}) = 0 .
\]
Taking into account (\ref{45}),  we calculate the approximate solution at the new time level. 

\textit{Stage 3}. Correction of the conditions at the boundaries of subdomains via the implicit scheme:
\[
   \frac{y^{n+1} - y^{n}}{\tau }
  + A  y^{n+1} = \varphi^n ,
  \quad \mathbf{x} \in \widehat{\omega} ,
\]
\[
  y^{n+1} = y^{n+1/2} ,
  \quad \mathbf{x} \notin  \widehat{\omega} .
\]

This numerical implementation (stages 1--3) of regionally-additive 
scheme  (\ref{45}), (\ref{48}), (\ref{49})
is nothing but the scheme of the domain decomposition 
\cite{1013.65106,1120.65098,1076.65079,1099.65084,1094.65097,1168.65387})
with the explicit-implicit procedure for calculating the boundary conditions at the boundaries of subdomains. 

Regionally-additive scheme (\ref{45}), (\ref{48}), (\ref{49}) 
has the first order approximation in  $\tau$. 
It is possible to use the schemes of second order where
\begin{equation}\label{52}
  \frac{y^{n+1/2} - y^{n}}{\tau/2 }
  + \chi_{1} A y^{n+1/2} + \chi_{2} A  y^{n} = \varphi^n,
\end{equation}
\begin{equation}\label{53}
  \frac{y^{n+1} - y^{n+1/2}}{\tau/2 }
  + \chi_{1} A  y^{n+1/2} + \chi_{2} A  y^{n+1} = \varphi^n .
\end{equation}
with $\varphi^n = f(t^{n+1/2})$.
Schemes  (\ref{48}), (\ref{49}) and (\ref{52}), (\ref{53})  we consider as the 
operator analogs of the classical schemes of alternating directions. 
They are special cases of more general factorized schemes. 

Consider the factorized scheme 
\begin{equation}\label{54}
  B_1 B_2 \frac{y^{n+1} - y^{n}}{\tau}
  + A  y^{n} = \varphi^n ,
\end{equation}
where 
\begin{equation}\label{55}
  B_{\alpha} = E + \sigma \tau \chi_{\alpha} A,
  \quad \alpha = 1, 2,
\end{equation} 
with the right-hand side specified in the form $\varphi^n = f(\sigma t^{n+1}+ (1-\sigma) t^{n})$.
Direct substitutions verify that scheme (\ref{54}), (\ref{55}) 
coincides with scheme  (\ref{48}), (\ref{49}) 
at $\sigma = 1$ and with scheme (\ref{52}), (\ref{53}) at $\sigma = 1/2$.

For the factorized scheme (\ref{54}), (\ref{55}) it is possible
to use the three-stage computational implementation with  explicit-implicit 
calculations of interface boundary conditions. 
We introduce, for example, the new grid function  $\widetilde{y}^{n+1}$
and instead of(\ref{54})  in view of (\ref{55}) we solve two differential equations: 
\begin{equation}\label{56}
  (E + \sigma \tau \chi_{1} A) \frac{\widetilde{y}^{n+1} - y^{n}}{\tau}
  + A  y^{n} = \varphi^n ,
\end{equation}
\begin{equation}\label{57}
  (E + \sigma \tau \chi_{2} A) \frac{y^{n+1} - y^{n}}{\tau} =
  \frac{\widetilde{y}^{n+1} - y^{n}}{\tau} .
\end{equation}
Taking into account (\ref{45}), we obtain from (\ref{56}) that
for nodes at common boundaries (\textit{Stage 1} --- the explicit scheme for boundary nodes):
\begin{equation}\label{58}
  \frac{\widetilde{y}^{n+1} - y^{n}}{\tau}
  + A  y^{n} = \varphi^n ,
  \quad \mathbf{x} \in \widehat{\omega} .
\end{equation}
For subdomains we have: 
\[
  (E + \sigma \tau A) \frac{y^{n+1} - y^{n}}{\tau} 
  + A  y^{n} = \varphi^n , 
  \quad \mathbf{x} \notin  \widehat{\omega} .
\]
This corresponds to (\textit{Stage 2} --- the implicit scheme in the subdomains) 
usage of the implicit scheme with weight  $\sigma$ for the difference solution in the subdomains. 
The implementation of (\ref{57}) (\textit{Stage 3} --- the implicit scheme for the boundary nodes)
in view of  (\ref{45}) is: 
\begin{equation}\label{59}
  (E + \sigma \tau A) \frac{y^{n+1} - y^{n}}{\tau} =
  \frac{\widetilde{y}^{n+1} - y^{n}}{\tau} ,
  \quad \mathbf{x} \in \widehat{\omega} ,  
\end{equation}
\[
  y^{n+1} = \widetilde{y}^{n+1} ,
  \quad \mathbf{x} \notin  \widehat{\omega} .
\]
Taking into account (\ref{58}), equation  (\ref{59}) can be written as 
\[
  (E + \sigma \tau A) \frac{y^{n+1} - y^{n}}{\tau} 
  + A  y^{n} = \varphi^n ,
  \quad \mathbf{x} \in \widehat{\omega} .
\]
At this stage all computational work is associated only with correction of 
the internal boundary conditions via this implicit scheme with weights. 

\begin{thm} 
\label{t-3} 
Factorized regionally-additive difference scheme (\ref{39}), (\ref{54}), (\ref{55})
is unconditionally stable for  $\sigma \geq 1/2$,
and for the difference solution the following estimate  holds
\begin{equation}\label{60}
  \|B_2 y^{n+1}\|_{A} \leq \|B_2 y^{n}\|_{A} + 
  \tau \|B_1^{-1}\varphi^n\|_A,
  \quad n = 0,1, ..., N-1 .
\end{equation}
\end{thm} 

\begin{pf} 
It is convenient firstly to symmetrize the factorized scheme (\ref{54}), (\ref{55}). 
Let $v^{n} = A^{1/2} y^{n}$ and
\[
  \widetilde{B}_{\alpha} = E + \sigma \tau \widetilde{A}_{\alpha},
  \quad  \widetilde{A}_{\alpha} = A^{1/2} \chi_{\alpha} A^{1/2},
  \quad \alpha = 1, 2 .
\]
Then equation (\ref{54}) can be rewritten as 
\begin{equation}\label{61}
  \widetilde{B}_1 \widetilde{B}_2 \frac{v^{n+1} - v^{n}}{\tau}
  + A  v^{n} = A^{1/2}\varphi^n  .  
\end{equation}
Assuming that $\widetilde{B}_2 v^{n} = w^{n}$,  from (\ref{61})
we obtain 
\begin{equation}\label{62}
  w^{n+1} = S w^{n} + \tau \widetilde{B}_1^{-1} A^{1/2}\varphi^n  , 
\end{equation}
where the operator of the transition to the new time level 
\begin{equation}\label{63}
  S = E - \tau \widetilde{B}_1^{-1} A \widetilde{B}_2^{-1} .
\end{equation}
Taking into account the above notation, from (\ref{63}) we obtain 
\[
  S = \frac{2\sigma - 1}{2\sigma} E + \frac{1}{2\sigma}
  \widetilde{B}_1^{-1} (\widetilde{B}_1 \widetilde{B}_2 - 
  2\sigma (\widetilde{A}_{1} + \widetilde{A}_{2}) ) \widetilde{B}_1^{-2}
\]
\[
  = \frac{2\sigma - 1}{2\sigma} E + S_1 S_2 ,
\]
where 
\[
  S_{\alpha} = (E + \sigma \tau \widetilde{A}_{\alpha})^{-1}
  (E - \sigma \tau \widetilde{A}_{\alpha}), 
  \quad \alpha = 1, 2 .
\]
If $\sigma \geq 0$, taking into account the non-negativity of the operators 
$\widetilde{A}_{\alpha},  \ \alpha = 1, 2$, we have 
\[
  \|S_{\alpha}\| \leq 1,
  \quad \alpha = 1, 2 .
\]
With stronger restrictions $\sigma \geq 1/2$ we find that 
$\|S\| \leq 1$. 
From (\ref{52}) we obtain the estimate 
\[
  \|w^{n+1}\| = \|w^{n}\| + \tau \|\widetilde{B}_1^{-1} A^{1/2}\varphi^n \|.  
\]
This is the required estimate (\ref{60}).
\end{pf} 
 
The fundamental issue in the construction of domain decomposition schemes 
for unsteady problems is to estimate the  convergence rate 
for the approximate solution. Accuracy depends on a computational 
grid (the width of the overlapping) and therefore regionally-additive schemes 
belong to the class of conditionally convergent. 
The situation can be illustrated by the example of the above 
factorized decomposition schemes (\ref{54}), (\ref{55}).

Analysis of the accuracy will be conducted in the standard way by considering 
the corresponding problem for the error 
\[
  z^{n}(\mathbf{x}) = y^{n}(\mathbf{x}) - u^n(\mathbf{x}),
  \quad \mathbf{x} \in \omega ,
\] 
where $u^n(\mathbf{x}) = u(\mathbf{x}, t^n)$ is  the exact solution of the differential problem 
(\ref{9})--(\ref{11}). 
From (\ref{39}), (\ref{54}), (\ref{55})  we obtain the problem for the error 
\begin{equation}\label{64}
  B_1 B_2 \frac{z^{n+1} - z^{n}}{\tau}
  + A  z^{n} = \psi^n ,
\end{equation}
\begin{equation}\label{65}
  z^0 = 0 .
\end{equation}
In view of (\ref{60}) for problem (\ref{64}), (\ref{65}) we have
\begin{equation}\label{66}
  \|B_2 z^{n+1}\|_{A} \leq \sum_{k=0}^{n}
  \tau \|B_1^{-1}\psi^k\|_A,
  \quad n = 0,1, ..., N-1 .
\end{equation}

For the approximation error we have 
\begin{equation}\label{67}
  \psi^n = \varphi^n - 
  B_1 B_2 \frac{u^{n+1} - u^{n}}{\tau}
  - A  u^{n} .
\end{equation}
Taking into account  (\ref{55}), from (\ref{67})  we obtain 
\[
  \psi^n =  \psi_1^n + \psi_2^n,
\]
\[
  \psi_1^n = \varphi^n -
  \left (E + \left (\sigma - \frac{1}{2} \right ) \tau A \right )
  \frac{u^{n+1} - u^{n}}{\tau}
  - A \frac{u^{n+1} + u^{n}}{2},
\]
\[
  \psi_2^n = - 
  \sigma^2 \tau^2 \chi_1 A \chi_2 A
  \frac{u^{n+1} - u^{n}}{\tau} .
\]
The first term of the error is the standard one for the schemes with weights, 
whereas the second term results from the splitting of subdomains. 
For sufficiently smooth solutions of problem (\ref{9})--(\ref{11}) we have 
\[
  \psi_1^n = \mathcal{O} (h^2 + \tau^2 +
  \left (\sigma - \frac{1}{2} \right ) \tau ) .  
\]
Let us consider the term $\psi_2^n$ in more detail.

Taking into account (\ref{66}) and introduced in the proof of Theorem~\ref{t-3} notation, we have 
\[
  \|B_1^{-1}\psi^n_2\|_A = 
  \|\widetilde{B}_1^{-1} A^{1/2}\psi^n \| 
\]
\[
  = \sigma^2 \tau^2 \left \| \widetilde{B}_1^{-1} A^{1/2} 
  \chi_1 A \chi_2 A   \frac{u^{n+1} - u^{n}}{\tau} \right \| 
\]
\[
  = \sigma \tau \left \| Q A^{1/2}\chi_2 A   \frac{u^{n+1} - u^{n}}{\tau} \right \| 
 \leq  \sigma \tau \left \| A^{1/2}\chi_2 A   \frac{u^{n+1} - u^{n}}{\tau} \right \| ,
\]
where 
\[
  Q = (E + \sigma \tau \widetilde{A}_{1})^{-1} 
  \sigma \tau \widetilde{A}_{1} .
\]
Thus 
\[
  \|B_1^{-1}\psi^n_2\|_A = 
  \mathcal{O} (\sigma \tau \|\chi_2\|_A) .
\]
These arguments allow us to formulate the following statement. 

\begin{thm} 
\label{t-4} 
For the error of the factorized regionally-additive difference scheme (\ref{39}), (\ref{54}), (\ref{55})
with $\sigma \geq 1/2$ we have for problem  (\ref{9})--(\ref{11}) 
the following estimate 
\begin{equation}\label{68}
  \|B_2 z^{n+1}\|_{A} \leq M \left ( h^2 + \tau^2 +
  \left (\sigma - \frac{1}{2} \right ) \tau 
  + \sigma \tau \|\chi_2\|_A \right ). 
\end{equation}
\end{thm} 

For considered here substructuring domain 
decomposition schemes with the grid  elliptic operators of second order (\ref{22})
and splitting (\ref{39}), (\ref{45}) estimate  (\ref{68}) gives 
\begin{equation}\label{69}
  \|B_2 z^{n+1}\|_{A} \leq M \left ( h^2 + \tau^2 +
  \left (\sigma - \frac{1}{2} \right ) \tau 
  + \sigma \tau \widehat{h}^{-1/2} h^{-1/2} \right ) .
\end{equation}
Note also that the use of the scheme with $\sigma = 1/2$ does not increase the order of accuracy. 
But in this case the main error term is two times lower compared to $\sigma = 1$.

A slightly different algorithm can be implemented. 
Instead of (\ref{54}) we apply the factorized scheme 
\begin{equation}\label{70}
  B_2 B_1 \frac{y^{n+1} - y^{n}}{\tau}
  + A  y^{n} = \varphi^n ,
\end{equation}
i.e. we comutate operators  $B_1$  and $B_2$.

The implementation of scheme (\ref{70})  will differ slightly from 
the implementation of  scheme (\ref{54}). 
Similar to (\ref{56}), (\ref{57}) we have 
\begin{equation}\label{71}
  (E + \sigma \tau \chi_{2} A) \frac{\widetilde{y}^{n+1} - y^{n}}{\tau}
  + A  y^{n} = \varphi^n ,
\end{equation}
\begin{equation}\label{72}
  (E + \sigma \tau \chi_{1} A) \frac{y^{n+1} - y^{n}}{\tau} =
  \frac{\widetilde{y}^{n+1} - y^{n}}{\tau} .
\end{equation}
At stage (\ref{71}) we use the implicit scheme for the nodes at the boundaries of 
the subdomains and explicit scheme in the subdomains. 
Note that for the explicit scheme it is enough to evaluate only the boundary nodes. 
At stage (\ref{72})  the solution in the subdomains is calculated using the implicit scheme. 
Thus, the computational cost in case of the factorized scheme (\ref{70}) 
remains practically the same as for scheme  (\ref{54}).

\section{Schemes of multi-component splitting}
\label{sec:6}

Constructed above factorized schemes of  the two-component splitting
can be generalized in various directions. 
The most fundamental issue is to construct such schemes in the case of general multi-component splitting. 
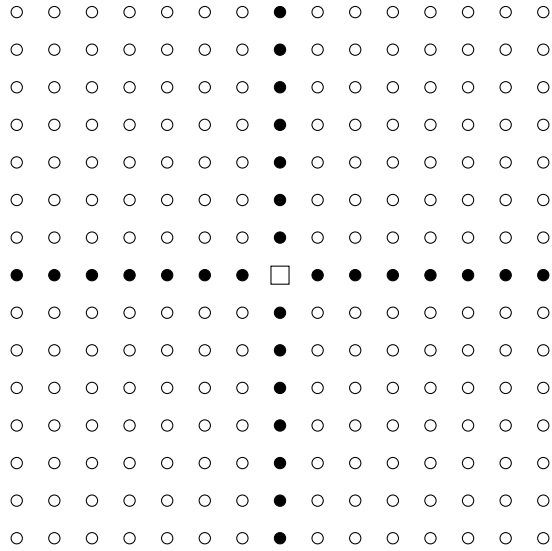
\begin{figure}[ht] 
  \begin{center}
	\begin{tikzpicture}[domain=0:4]
		\foreach \x in {1,...,15}
  			\foreach \y in {1,...,15}
	  			\draw [{fill=white}] (0.5*\x,0.5*\y) circle (0.075);
		\foreach \x in {1,...,15}
	  		\fill [black] (0.5*\x,4) circle (0.075);
		\foreach \y in {1,...,15}
	  		\fill [black] (4,0.5*\y) circle (0.075);
	  	\draw [{fill=white}] (4-0.12,4-0.12) rectangle +(0.24,0.24);
	\end{tikzpicture}
    \caption{Three-component decomposition without the overlapping of subdomains} 
	\label{f-3}
  \end{center}
\end{figure}

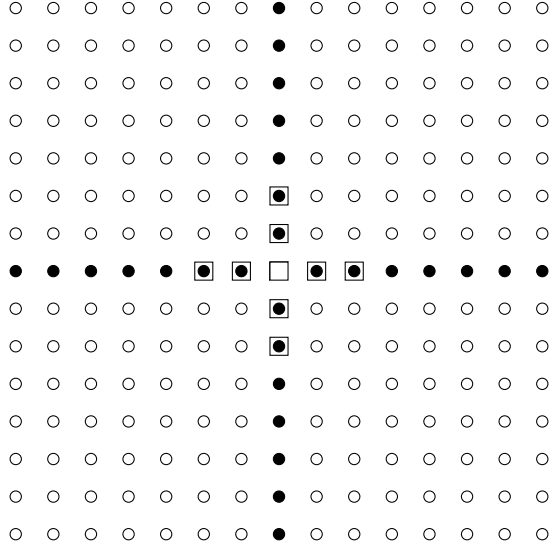
\begin{figure}[ht] 
  \begin{center}
	\begin{tikzpicture}[domain=0:4]
		\foreach \x in {1,...,15}
  			\foreach \y in {1,...,15}
	  			\draw [{fill=white}] (0.5*\x,0.5*\y) circle (0.075);
		\foreach \x in {1,...,15}
	  		\fill [black] (0.5*\x,4) circle (0.075);
		\foreach \y in {1,...,15}
	  		\fill [black] (4,0.5*\y) circle (0.075);
		\foreach \x in {0,...,4}
	  		\draw (0.5*\x+3-0.12,4-0.12) rectangle +(0.24,0.24);
		\foreach \y in {0,...,4}
	  		\draw (4-0.12,0.5*\y+3-0.12) rectangle +(0.24,0.24);
	  	\draw [{fill=white}] (4-0.12,4-0.12) rectangle +(0.24,0.24);
	\end{tikzpicture}
    \caption{Three-component decomposition with the overlapping of subdomains } 
	\label{f-4}
  \end{center}
\end{figure}

The need for such an extension results from, in particular, 
calculations of conditions at the boundaries of subdomains, i.e. 
the solution of problems on graphs for two-dimensional problems. 
In the considered two-dimensional problems in a rectangle and rectangular grids,
the implementation of, for example, (\ref{49})  does not face significant problems. 
However, for more general situations, for example, for three-dimensional boundary 
value problems, the solution of these grid problems can be difficult.  
Such considerations lead us to the need of constructing procedures of
decomposition for the set of boundary nodes of subdomains. 
A characteristic example is shown in Fig.~\ref{f-3}. 
The set of boundary nodes is divided into two parts:  
$\widehat{\omega} = \widehat{\omega}_{s}  \cup \widehat{\omega}_{m}$.
Here the set of nodes at the boundary of  two subdomains is denoted as
$\widehat{\omega}_{s}$ (in Fig.~\ref{f-3} it is depicted as  $\bullet$).
The set of nodes that lie at the boundaries of a greater number of subdomains is 
designated as $\widehat{\omega}_{m}$ (in Fig.~\ref{f-3} it is presented as  {\tiny $\square$}).

Instead of the two-component splitting (\ref{39}), (\ref{45}), 
we use now the three-component  splitting  (\ref{39}) with $ p = 3 $ and
\[
  \chi_{2}(\mathbf{x}) = \left \{
  \begin{array}{ll}
    1, &  \mathbf{x} \in \widehat{\omega}_{s}, \\
    0, &  \mathbf{x} \notin  \widehat{\omega}_{s}, \\
  \end{array}
  \right .
  \quad 
  \chi_{3}(\mathbf{x}) = \left \{
  \begin{array}{ll}
    1, &  \mathbf{x} \in \widehat{\omega}_{m}, \\
    0, &  \mathbf{x} \notin  \widehat{\omega}_{m}, \\
  \end{array}
  \right .
\]
\begin{equation}\label{73}
  \chi_{1}(\mathbf{x}) = 1 - \chi_{2}(\mathbf{x}) - \chi_{3}(\mathbf{x}),
  \quad \mathbf{x} \in \omega .
\end{equation}
With such a decomposition calculations in some parts of the  subdomain boundaries
(on the set  $\widehat{\omega}_{s}$) can be performed independently 
using known conditions at the nodes of crossing (on the set  $\widehat{\omega}_{m}$).

Local computations of of the solution at boundary crossings introduces additional errors. 
To improve the accuracy of the approximate solution at the boundaries of 
subdomains, it is possible to apply algorithms with the overlapping of subdomains. 
Such a situation at  the grid level is shown in Fig.~\ref{f-4}.
There is highlighted the set of boundary nodes  $\widehat{\omega}_{m}$, 
which lie near the boundary crossing  and $\widehat{\omega}_{s}  \cap \widehat{\omega}_{m} \neq \varnothing$.
With this in mind, instead of (\ref{73}) we set 
\[
  \chi_{2}(\mathbf{x}) = \left \{
  \begin{array}{ll}
    > 0, &  \mathbf{x} \in \widehat{\omega}_{s}, \\
    0, &  \mathbf{x} \notin  \widehat{\omega}_{s}, \\
  \end{array}
  \right .
  \quad 
  \chi_{3}(\mathbf{x}) = \left \{
  \begin{array}{ll}
    > 0, &  \mathbf{x} \in \widehat{\omega}_{m}, \\
    0, &  \mathbf{x} \notin  \widehat{\omega}_{m}, \\
  \end{array}
  \right .
\]
\begin{equation}\label{74}
  \chi_{2}(\mathbf{x}) + \chi_{3}(\mathbf{x}) = 1,
  \quad \mathbf{x} \in \widehat{\omega} ,
  \quad \chi_{1}(\mathbf{x}) = 1 - \chi_{2}(\mathbf{x}) - \chi_{3}(\mathbf{x}),
  \quad \mathbf{x} \in \omega .
\end{equation}

For the general multi-component ($p > 2$) decomposition it is possible to construct in a more simple way 
regularized additive schemes \cite{0963.65091,samarskii1998regularized}.
For solving problem  (\ref{25}), (\ref{40}), (\ref{42}) 
we can use the additive scheme of full approximation 
\begin{equation}\label{75}
  \frac{y^{n+1} - y^{n}}{\tau }
  + \widetilde{A} y^{n} 
  = \varphi^n ,
\end{equation}
where 
\begin{equation}\label{76}
  \widetilde{A} = \sum_{\alpha =1}^{p} \widetilde{A}_{\alpha},
  \quad \widetilde{A}_{\alpha} = 
  (E + \sigma \tau \chi_{\alpha} A)^{-1} \chi_{\alpha} A ,
  \quad \alpha =1,2, ...,p .
\end{equation}
This scheme is characterized by the fact that each operator term $\chi_{\alpha} A, \ \alpha =1,2, ...,p$
is perturbed with an error $\mathcal{O}(\tau)$.

\begin{thm} 
\label{t-5} 
Regularized difference scheme (\ref{75}), (\ref{76}) is unconditionally stable for 
 $\sigma \geq p/2$,
and for the difference solution we have the estimate 
\begin{equation}\label{77}
  \|y^{n+1}\|_{A} \leq \|y^{n}\|_{A} + 
  \tau \|\varphi^n\|_{A},
  \quad n = 0,1, ..., N-1 .
\end{equation}
\end{thm} 

\begin{pf} 
The operator $\widetilde{A}$ can be written in the form 
\begin{equation}\label{78}
  \widetilde{A} = \sum_{\alpha =1}^{p} A^{-1/2} R_{\alpha} A^{1/2},
\end{equation}
where 
\[
  R_{\alpha} = (E + \sigma \tau C_{\alpha})^{-1} C_{\alpha},
  \quad C_{\alpha} = A^{1/2} \chi_{\alpha} A^{1/2},
  \quad \alpha =1,2, ...,p .
\]
Thus 
\[
  C_{\alpha} = C_{\alpha}^* \geq 0,
  \quad R_{\alpha} = R_{\alpha}^* \geq 0,
  \quad \alpha =1,2, ...,p .
\]
With this in mind the difference scheme  (\ref{75}), (\ref{76})
is written as 
\begin{equation}\label{79}
  \frac{v^{n+1} - v^{n}}{\tau }
  + \sum_{\alpha =1}^{p} R_{\alpha} v ^{n} 
  = A^{1/2} \varphi^n ,
\end{equation}
where, as before, $v^n = A^{1/2} y^n$.
From (\ref{79}) we have 
\begin{equation}\label{80}
  v^{n+1} = S v^{n} + \tau A^{1/2} \varphi^n   
\end{equation}
with the transition operator 
\[
  S = E  - \tau \sum_{\alpha =1}^{p} R_{\alpha} .
\]
Using this representation, we set 
\begin{equation}\label{81}
  S = \frac{1}{p}\sum_{\alpha =1}^{p} S_{\alpha} ,
  \quad S_{\alpha} = E  - p \tau R_{\alpha} ,
  \quad \alpha =1,2, ...,p .
\end{equation}
For the individual terms with  the above notation we obtain 
\[
  S_{\alpha} = (E + \sigma \tau C_{\alpha})^{-1}
  (E + (\sigma - p) \tau C_{\alpha}).
\]
Under the constraints  $\sigma \geq p/2$  we have 
$\|S_{\alpha} \| \leq 1$, which allow us to obtain from  (\ref{80}) the following estimate 
\[
  \|v^{n+1}\| \leq \|v^{n} \| + \tau \|A^{1/2} \varphi^n \|,
\] 
which is nothing but  (\ref{77}).
\end{pf} 

Standard finite-difference schemes of component-wise splitting 
\cite{0971.65076,Marchuk:1990:SAD,0209.47103} can be easy constructed using
the transition operator. 
With regard to our problem, we shall again start with notation (\ref{80}), 
but instead the additive structure (see (\ref{81}) employ  multiplicative one 
\begin{equation}\label{82}
  S = \prod_{\alpha =p}^{1} S_{\alpha} ,
  \quad S_{\alpha} = E  -  \tau R_{\alpha} ,
  \quad \alpha =1,2, ...,p .
\end{equation}
In this case we have $\|S_{\alpha} \| \leq 1$ for $\sigma \geq 1/2$. 

The implementation of  component-wise splitting scheme is performed as
a sequence of intermediate difference problems similar to (\ref{80}): 
\begin{equation}\label{83}
  v^{n+\alpha/p} = S_{\alpha} v^{n+(\alpha-1)/p} + 
  \tau A^{1/2} \varphi_{\alpha}^n ,  
  \quad \alpha =1,2, ...,p .
\end{equation}
Comparing with (\ref{80}), (\ref{81}), we obtain 
\[
  \varphi^n = \sum_{\alpha =p}^{1} 
  \prod_{\beta =p}^{\alpha -1} S_{\beta} \varphi_{\alpha}^n .
\]
Without loss of accuracy, we can consider only the simplest choice for 
$\varphi_{\alpha}^n, \ \alpha =1,2, ...,p$, where
\begin{equation}\label{84}
  \varphi^n = \sum_{\alpha =p}^{1} \varphi_{\alpha}^n .
\end{equation}

With this notation the difference equations (\ref{83}) can be written as follows
\[
  \frac{y^{n+\alpha/p} - y^{n+(\alpha-1)/p}}{\tau}
  + (E + \sigma \tau \chi_{\alpha} A)^{-1} \chi_{\alpha} A y^{n+(\alpha-1)/p}
  = \varphi_{\alpha}^n ,
\]
\begin{equation}\label{85}
  \alpha =1,2, ...,p .
\end{equation}
We can formulate now the following statement. 

\begin{thm} 
\label{t-6} 
Additive component-wise splitting scheme (\ref{84}), (\ref{85})
is unconditionally stable for $\sigma \geq 1/2$,
and for the difference solution estimate (\ref{77}) holds.
\end{thm} 

The considering scheme has the first-order approximation in time. 
However, in the case of two-component splitting at $\sigma = 1/2$ the approximation error is $\mathcal{O}(\tau^2)$
(See, e.g., \cite{samgul}). 
This variant is used in  \cite{dryja2007} to construct domain decomposition schemes. 
For the general multi-component splitting the additive schemes of second order in time 
are based on the symmetrization of transition operator \cite{0184.38503,fryazinov1968economical}. 
In this case, instead of (\ref{82}) we can use, for example, 
\[
  S = \prod_{\beta  =1}^{p} S_{\beta} \prod_{\alpha =p}^{1} S_{\alpha} ,
  \quad   S_{\alpha} = \left (E + \frac{\tau}{4} C_{\alpha} \right )^{-1}
  \left (E - \frac{\tau}{4} C_{\alpha} \right ),
  \quad \alpha =1,2, ...,p .
\]
Thus we make two half-steps in time for  $\sigma \geq 1/2$
in sequential solving problems for operators 
$\chi_{\alpha} A, \ \alpha =1,2, ...,p$, and then for the operators $\chi_{\beta} A, \ \beta =p,p-1, ...,1$.

Regularized scheme (\ref{75}), (\ref{76})  can be written in the form similar to  (\ref{85}):
\[
  \frac{y^{n+\alpha/p} - y^{n+(\alpha-1)/p}}{\tau}
  + (E + \sigma \tau \chi_{\alpha} A)^{-1} \chi_{\alpha} A y^{n}
  = \varphi_{\alpha}^n ,
\]
\begin{equation}\label{86}
  \alpha =1,2, ...,p .
\end{equation}
In contrast to  (\ref{85}) here  the obtained $y^{n+(\alpha-1)/p}$ is used  
for solving the problem for $y^{n+\alpha/p}$ only partially. 
This increasing of explicitness results in 
a more strong condition of stability 
(instead of $\sigma \geq 1/2$  we have  $\sigma \geq p/2$).

The numerical implementation of  scheme (\ref{84}), (\ref{85})
is shown below for decomposition (\ref{45}). 
In accordance with (\ref{84}) we set 
\[
  \varphi_{\alpha}^n = \chi_{\alpha} \varphi^n ,
  \quad \alpha =1,2, ...,p .
\]
From (\ref{85})  we obtain 
\begin{equation}\label{87}
  (E + \sigma \tau \chi_{1} A) 
  \frac{y^{n+1/2} - y^{n}}{\tau} + \chi_{1} A y^{n}
  = (E + \sigma \tau \chi_{1} A) \chi_{\alpha} \varphi^n ,
\end{equation}
\begin{equation}\label{88}
  (E + \sigma \tau \chi_{2} A) 
  \frac{y^{n+1} - y^{n+1/2}}{\tau} + \chi_{2} A y^{n+1/2}
  = (E + \sigma \tau \chi_{2} A) \chi_{2} \varphi^n .
\end{equation} 
In finding  $y^{n+1/2}$  from (\ref{87}) (\textit{Stage 1}) 
we solve boundary value problems in subdomains using the implicit scheme. 
The boundary conditions are taken from the previous time level, i.e., 
\[
  (E + \sigma \tau A) \frac{y^{n+1/2} - y^{n}}{\tau }
  + A  y^{n} = (E + \sigma \tau A) \varphi^n ,
  \quad \mathbf{x} \notin  \widehat{\omega} .
\]
\[
  y^{n+1/2} = y^{n} ,
  \quad \mathbf{x} \in \widehat{\omega} .
\]
Conditions at the common boundaries are corrected during evaluation  $y^{n+1}$ from (\ref{88}) (\textit{Stage 2}):
\[
  (E + \sigma \tau A) \frac{y^{n+1} - y^{n+1/2}}{\tau }
  + A  y^{n+1/2} = (E + \sigma \tau A) \varphi^n ,
  \quad \mathbf{x} \in \widehat{\omega} ,
\]
\[
  y^{n+1} = y^{n+1/2} ,
  \quad \mathbf{x} \notin  \widehat{\omega} .
\]
The numerical implementation of the component-wise splitting scheme (\ref{84}), (\ref{85}) 
is slightly reduced in compare with the factorized domain decomposition scheme 
(\ref{48}), (\ref{49}) (there is no explicit calculations of the interface boundary conditions). 
Similarly, two-stage implementation takes place for regularized scheme (\ref{84}), (\ref{86}).
 
\section{Hyperbolic equations of second order}
\label{sec:7}

Possibilities of constructing domain decomposition schemes to solve
boundary value problems  for hyperbolic equation  of second order (\ref{10}),
(\ref{17}), (\ref{18})  are more restricted. 
Here we note only the regularized schemes similar to (\ref{75}), (\ref{76})
for the parabolic problem (\ref{9})--(\ref{11}).
\begin{equation}\label{89}
  \frac{y^{n+1} - 2y^{n} + y^{n-1}}{\tau^2}
  + \widetilde{A} y^{n} 
  = \varphi^n ,
\end{equation}
where
\begin{equation}\label{90}
  \widetilde{A} = \sum_{\alpha =1}^{p} \widetilde{A}_{\alpha},
  \quad \widetilde{A}_{\alpha} = 
  (E + \sigma \tau^2 \chi_{\alpha} A)^{-1} \chi_{\alpha} A ,
  \quad \alpha =1,2, ...,p .
\end{equation}
This scheme has the second order of accuracy in time. 
The following statement is true. 

\begin{thm} 
\label{t-7} 
Regularized difference scheme (\ref{34}), (\ref{89}), (\ref{90})
is unconditionally stable for  $\sigma \geq p/4$,
and for the difference solution the following estimate is satisfied
\begin{equation}\label{91}
  S^{n+1} \leq \exp(\tau) S^{n} 
  + \frac{\tau^2}{2} \frac{\exp(\tau)}{\exp(0.5\tau)-1} \|\varphi^n \|^2_{D^{-1}},
  \quad n = 1,2, ..., N-1 ,
\end{equation}
where 
\[
  S^n =  \left \|\frac{y^n - y^{n-1}}{\tau } \right \|^2_D + 
  \left \|\frac{y^n + y^{n-1}}{2} \right \|_{A \widetilde{A}}^2 ,
\]
\[
  D = D^* = A \left (E - \frac{\tau^2}{4} \widetilde{A} \right ) .
\]
\end{thm} 

\begin{pf} 
The proof is conducted similarly to Theorem~\ref{t-2}.
Similarly  (\ref{78}), for the operator $\widetilde{A}$, taking into account (\ref{90}) , we have representation 
\begin{equation}\label{92}
  \widetilde{A} = \sum_{\alpha =1}^{p} A^{-1/2} R_{\alpha} A^{1/2},
\end{equation}
where now 
\[
  R_{\alpha} = (E + \sigma \tau^2 C_{\alpha})^{-1} C_{\alpha},
  \quad C_{\alpha} = A^{1/2} \chi_{\alpha} A^{1/2},
  \quad \alpha =1,2, ...,p 
\]
with self-adjoint and non-negative operators $C_{\alpha}$ è $R_{\alpha}, \ \alpha =1,2, ...,p$.

Difference scheme (\ref{89}), (\ref{90}) can be written as 
\begin{equation}\label{93}
  \frac{v^{n+1} - 2v^{n} - v^{n-1}}{\tau }
  + \sum_{\alpha =1}^{p} R_{\alpha} v ^{n} 
  = A^{1/2} \varphi^n 
\end{equation}
for $v^n = A^{1/2} y^n$.
Using the notation 
\[
  \zeta^n = \frac{v^n + v^{n-1}}{2} ,
  \quad \eta^n = \frac{v^n - v^{n-1}}{\tau } ,
\]
write  (\ref{92})  as 
\begin{equation}\label{94}
  \left (E - \frac{\tau^2}{4}  R \right ) 
  \frac{\eta^{n+1} - \eta^{n}}{\tau} + 
  R \frac{\zeta^{n+1} + \zeta^{n}}{2} = A^{1/2} \varphi^n ,
\end{equation}
where 
\[
  R = R^* = \sum_{\alpha =1}^{p} R_{\alpha} .
\]
Multiply (\ref{94})  scalarly in  $H$ by 
\[
  2 (\zeta^{n+1} - \zeta^{n}) = \tau (\eta^{n+1} + \eta^{n}) 
\]
and obtain the equality 
\begin{equation}\label{95}
  S^{n+1} - S^{n} =
  \tau (A^{1/2} \varphi^n, (\eta^{n+1} + \eta^{n})) ,
\end{equation}
where 
\[
  S^n =  \|\eta^{n}\|^2_{\widetilde{D}} + \|\zeta^{n}\|_R^2 
  = \left \|\frac{v^n - v^{n-1}}{\tau } \right \|^2_{\widetilde{D}} + 
  \left \|\frac{v^n + v^{n-1}}{2} \right \|_R^2 ,
\]
\[
  \widetilde{D} = E - \frac{\tau^2}{4}  R .
\]
and $\widetilde{D} > 0$.
For the first term we have 
\[
   \left \|\frac{v^n - v^{n-1}}{\tau } \right \|^2_{\widetilde{D}} =
   \left \|\frac{y^n - y^{n-1}}{\tau } \right \|^2_{D} 
\]
that results from
\[
  A^{1/2} \widetilde{D} A^{1/2} = A - \frac{\tau^2}{4} A^{1/2} R A^{1/2}
  = A \left (E - \frac{\tau^2}{4} \widetilde{A} \right ) = D .
\]

To prove the inequality $\widetilde{D} > 0$  for $\sigma \geq p/4$, we set 
\[
  \widetilde{D} = \frac{1}{p} \sum_{\alpha =1}^{p} \widetilde{D}_{\alpha},
  \quad \widetilde{D}_{\alpha} = E - \frac{\tau^2}{4} p R_{\alpha},
  \quad \alpha = 1,2,...,p.
\]
For each individual term we have $\widetilde{D}_{\alpha}, \ \alpha = 1,2,...,p$, if
\[
  E - \frac{\tau^2}{4} p (E + \sigma \tau^2 C_{\alpha})^{-1} C_{\alpha} > 0.
\]
We have
\[
  E + \sigma \tau^2 C_{\alpha} - \frac{\tau^2}{4} p C_{\alpha} > E .
\]
for $\sigma \geq p/4$ for each  $\alpha = 1,2,...,p$.

For the right-hand side of  (\ref{95}) with $\widetilde{\varphi}^n = A^{1/2} \varphi^n$, we use estimates 
\[
  \tau (\widetilde{\varphi},(\eta^{n+1} + \eta^{n})) \leq 
  \frac{\tau  }{2 \varepsilon}    \|\eta^{n+1} + \eta^{n}\|^2_{\widetilde{D}} +
  \frac{\tau}{2}  \varepsilon \|\widetilde{\varphi}^n \|^2_{\widetilde{D}^{-1}} ,
\]
\[
  \|\eta^{n+1} + \eta^{n}\|^2_{\widetilde{D}} \leq 
  2 (\|\eta^{n+1}\|^2_{\widetilde{D}}  + \|\eta^{n}\|^2_{\widetilde{D}} ) .
\]
We obtain the inequality 
\begin{equation}\label{96}
  \left (1 - \frac{\tau }{\varepsilon } \right )   S^{n+1} \leq 
  \left (1 + \frac{\tau }{\varepsilon } \right ) S^{n} +
  \frac{\tau}{2}  \varepsilon  \|\widetilde{\varphi}^n \|^2_{\widetilde{D}^{-1}} .
\end{equation}
We assume that 
\[
  1 - \frac{\tau }{\varepsilon } = \exp(-0.5 \tau) , 
\] 
and therefore 
\[
  1 + \frac{\tau }{\varepsilon } < \exp(0.5 \tau) .
\]
With our notation it is easy to obtain from (\ref{96}) 
the required estimate of stability (\ref{91}).
\end{pf} 

This estimate of stability is characterized by using
more complex norms in compare with the case of standard schemes with weights 
(compare (\ref{35}) and (\ref{91})).
The numerical implementation of the regularized scheme (\ref{89}), (\ref{90})
can be performed similarly to scheme (\ref{75}), (\ref{76}).

\section{Numerical results for model problems}
\label{sec:8}

Numerical experiments are performed for the parabolic equation (\ref{9}), where 
\begin{equation}\label{97}
  k(\mathbf{x}) = 1,
  \quad f(\mathbf{x},t) = 0,
   \quad {\bf x}\in \Omega,
   \quad 0 < t \leq T .
\end{equation}
Problem (\ref{9})--(\ref{11}), (\ref{96})  is considered in the unit square  $l_1 = l_2 = 1$,
and the initial condition has the form 
\begin{equation}\label{98}
  u^0(\mathbf{x}) = 
  \sin(n_1 \pi x_1) \sin(n_2 \pi x_2) ,
   \quad {\bf x}\in \Omega ,
\end{equation}
for natural $n_1$  and $n_2$.
The solution of problem (\ref{9})--(\ref{11}), (\ref{97}), (\ref{98}) is written as 
\begin{equation}\label{99}
  u(\mathbf{x},t) = 
  \exp (- \pi^2(n_1^2  + n_2^2) t) 
  \sin(n_1 \pi x_1) \sin(n_2 \pi x_2) .
\end{equation}

\begin{figure}[h]
  \includegraphics[width=0.8\textwidth,angle=-0]{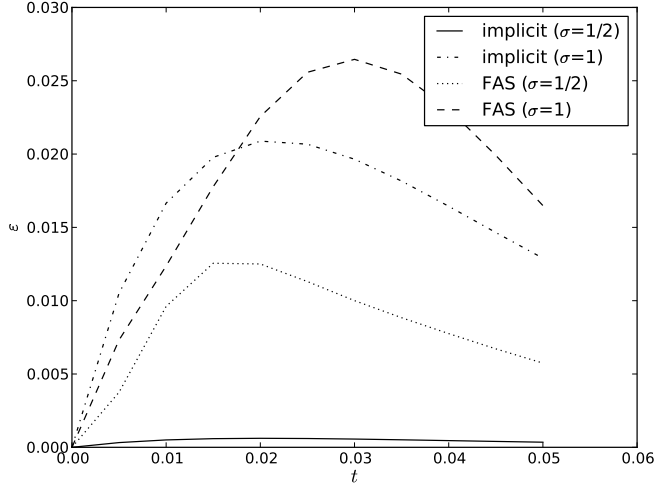}
  \caption{The error of factorized regionally-additive scheme}
\label{f-5}
\end{figure}

\begin{figure}[h]
  \includegraphics[width=0.8\textwidth,angle=-0]{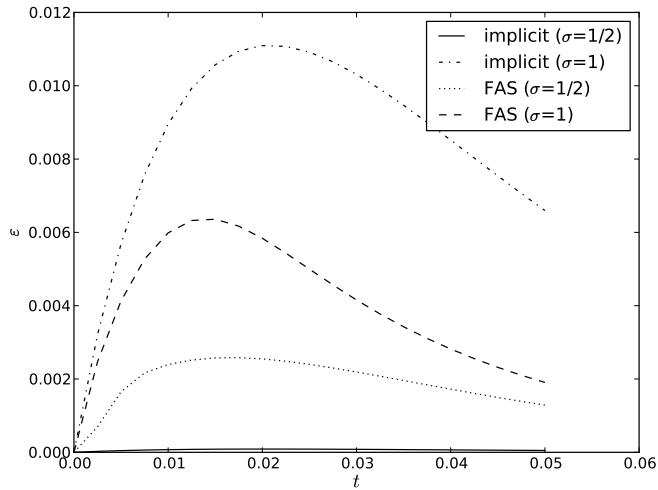}
  \caption{Reducing of the time step ($\tau = 0.005$)}
\label{f-6}
\end{figure}

\begin{figure}[h]
  \includegraphics[width=0.8\textwidth,angle=-0]{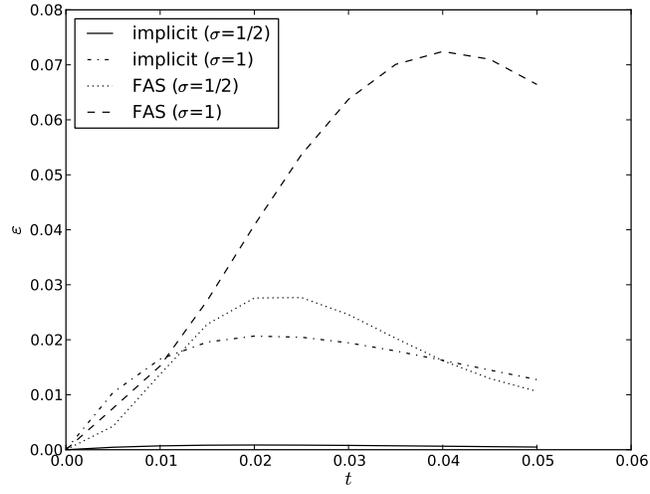}
  \caption{Reducing of the spatial step ($h = 1/80$)}
\label{f-7}
\end{figure}

\begin{figure}[h]
  \includegraphics[width=0.8\textwidth,angle=-0]{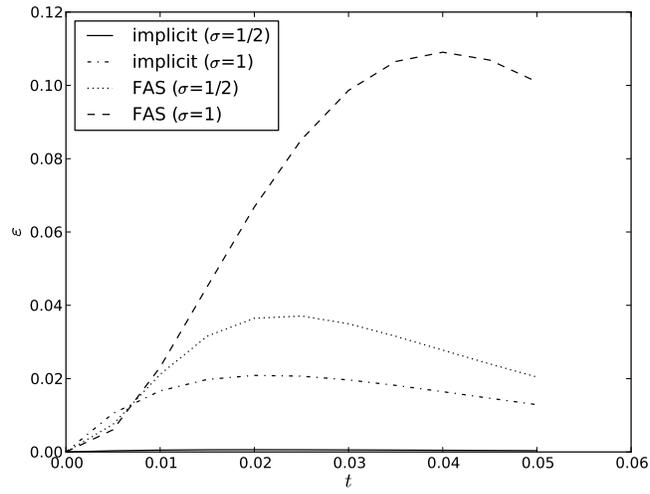}
  \caption{Increasing of the number of subdomains  ($\widehat{h} = 0.25$)}
\label{f-8}
\end{figure}

\begin{figure}[h]
  \includegraphics[width=0.8\textwidth,angle=-0]{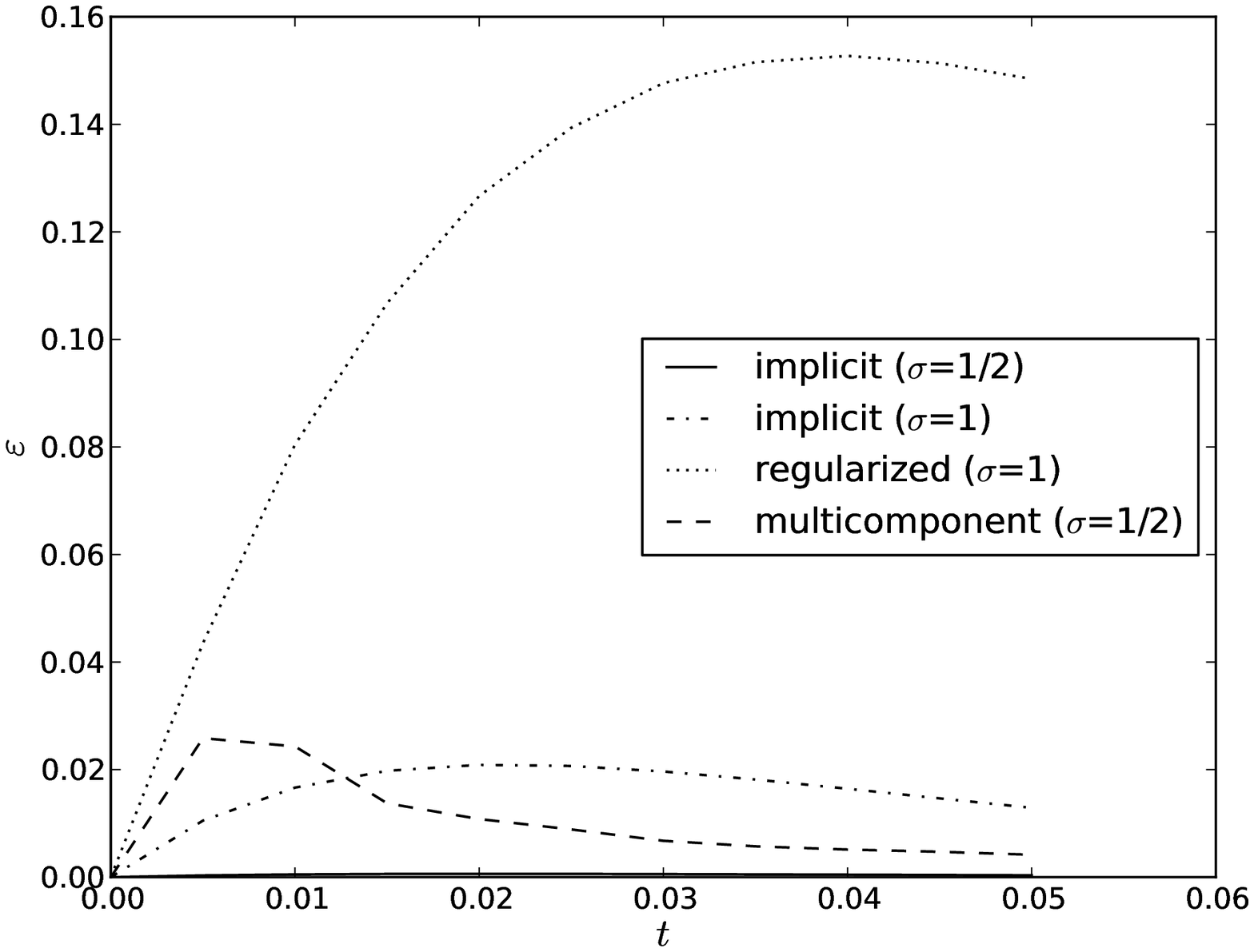}
  \caption{The error of the regularized and regionally-additive component-wise splitting schemes}
\label{f-9}
\end{figure}

The numerical  results derived using the regionally-additive schemes 
are compared with the difference solution obtained by means of 
the implicit scheme (\ref{30}), (\ref{31}) for   $\sigma = 1/2$ and $\sigma = 1$.
The error of the approximate solution was estimated as $\varepsilon(t^n) = \| y^n(\mathbf{x}) - u(\mathbf{x},t^n)\|$ 
at each particular time level. 

In the basic case we used $n_1=2, n_2=1$,
$N_1 = N_2$, $h_1 = h_2 = h = 1/40$, $T = 0.05, N = 10, \tau = 0.01$. 
The decomposition is carried out by cutting the  $\Omega$ into four squares 
($\widehat{h} = 0.5$).

Results of solving the test problem using the standard implicit 
schemes with weights (\ref{30}), (\ref{31}) of  second  ($\sigma = 1/2$) 
and first  ($\sigma = 1$)  orders of accuracy with respect to  $\tau$  
and factorized regionally-additive scheme (\ref{39}), (\ref{54}), (\ref{55}) 
(FAS) for the same values of the weight parameter $\sigma$  are presented in Fig.\ref{f-5}.

With the selected parameters the domain decomposition scheme,
constructed on the basis of classical factorized schemes, 
yields the approximate solution with a slightly larger error 
than the standard two-level scheme with weights. 
The effect of time step is presented in Fig. \ref{f-6},
where data are obtained with redusing time step ($\tau = 0.005$).
A more interesting effect is connected with the discretization in space (Fig.\ref{f-7}).
The effect of conditional convergence becomes more evident  for the regionally-additive 
scheme at $\sigma =1/2$.

In the study of the decomposition schemes particular attention 
should be paid to the dependence of accuracy of the approximate solutions on the number of subdomains.  
The error of the schemes for the increased number of subdomains (four times  ($\widehat{h} = 0.25$))
is shown in Fig.\ref{f-8}). 
Decreasing of the accuracy (compare Fig.\ref{f-5} and Fig.\ref{f-8}) is more significant
for the factorized regionally-additive schemes at $\sigma =1/2$.

For the model problem  (\ref{9})--(\ref{11}), (\ref{97}), (\ref{98})
with the exact solution (\ref{99}) and decomposition  (\ref{40}), (\ref{45})
we performed calculatuins via the above schemes of multicomponent-wise splitting. 
The error of the approximate solution for the basic case derived using the regularized 
regionally-additive scheme (\ref{75}), (\ref{76}) with  $\sigma =1$ 
is compared with the results of regionally-additive component-wise splitting scheme (\ref{84}), (\ref{85}) 
with $\sigma =1/2$  in Fig.\ref{f-9}.
The accuracy of the regularized scheme is clearly much lower. As for
the accuracy of the component-wise splitting  scheme, it is practically the same as 
the accuracy of the factorized regionally-additive scheme
(see Fig.\ref{f-5} and Fig.\ref{f-9}).

\clearpage

% \bibliographystyle{elsarticle-num}
% \bibliography{DDM-JCP}

\section*{References}

\begin{thebibliography}{10}
\expandafter\ifx\csname url\endcsname\relax
  \def\url#1{\texttt{#1}}\fi
\expandafter\ifx\csname urlprefix\endcsname\relax\def\urlprefix{URL }\fi
\expandafter\ifx\csname href\endcsname\relax
  \def\href#1#2{#2} \def\path#1{#1}\fi

\bibitem{pre05281749}
T.~Mathew, Domain decomposition methods for the numerical solution of partial
  differential equations, Lecture Notes in Computational Science and
  Engineering 61. Berlin: Springer. xiii, 764~p., 2008.

\bibitem{0931.65118}
A.~Quarteroni, A.~Valli, Domain decomposition methods for partial differential
  equations, Numerical Mathematics and Scientific Computation. Oxford:
  Clarendon Press. xv, 360 p., 1999.

\bibitem{0857.65126}
B.~Smith, Domain decomposition. Parallel multilevel methods for elliptic
  partial differential equations, Cambridge: Cambridge University Press. xii,
  224 p., 1996.

\bibitem{1069.65138}
A.~Toselli, O.~Widlund, Domain decomposition methods -- algorithms and theory,
  Springer Series in Computational Mathematics 34. Berlin: Springer. xv,
  450~p., 2005.

\bibitem{Cai:1991:ASA}
X.-C. Cai, Additive {S}chwarz algorithms for parabolic convection-diffusion
  equations, Numer. Math. 60~(1) (1991) 41--61.

\bibitem{Cai:1994:MSM}
X.-C. Cai, Multiplicative {S}chwarz methods for parabolic problems, SIAM J. Sci
  Comput. 15~(3) (1994) 587--603.

\bibitem{0825.65066}
Y.~A. Kuznetsov, New algorithms for approximate realization of implicit
  difference schemes, Sov. J. Numer. Anal. Math. Model. 3~(2) (1988) 99--114.

\bibitem{0766.65089}
Y.~A. Kuznetsov, Overlapping domain decomposition methods for fe-problems with
  elliptic singular perturbed operators, Fourth international symposium on
  domain decomposition methods for partial differential equations, Proc. Symp.,
  Moscow/Russ. 1990, 223-241 (1991) (1991).

\bibitem{1018.65103}
A.~A. Samarskii, P.~P. Matus, P.~N. Vabishchevich, Difference schemes with
  operator factors, Mathematics and its Applications (Dordrecht). 546.
  Dordrecht: Kluwer Academic Publishers. x, 384 p., 2002.

\bibitem{Laevsky}
Y.~M. Laevsky, Domain decomposition methods for the solution of two-dimensional
  parabolic equations, in: Variational-difference methods in problems of
  numerical analysis, no.~2, Comp. Cent. Sib. Branch, USSR Acad. Sci.,
  Novosibirsk, 1987, pp. 112--128, in Russian.

\bibitem{0719.65072}
P.~N. Vabishchevich, Difference schemes with domain decomposition for solving
  non-stationary problems, U.S.S.R. Comput. Math. Math. Phys. 29~(6) (1989)
  155--160.

\bibitem{0963.65091}
A.~A. Samarskii, P.~N. Vabishchevich, Additive schemes for problems of
  mathematical physics (Additivnye skhemy dlya zadach matematicheskoj fiziki),
  Moscow: Nauka. 320 p., 1999, in Russian.

\bibitem{0928.65102}
A.~A. Samarskii, P.~N. Vabishchevich, Factorized finite-difference schemes for
  the domain decomposition in convection-diffusion problems, Differ. Equations
  33~(7) (1997) 972--979.

\bibitem{vab_255}
A.~A. Samarskii, P.~N. Vabishchevich, Domain decomposition methods for
  parabolic problems, in: C.-H. Lai, P.~Bjorstad, M.~Gross, O.~Widlund (Eds.),
  Eleventh International Conference on Domain Decomposition Methods, DDM.org,
  1999, pp. 341--347.

\bibitem{0888.65097}
P.~N. Vabishchevich, Finite-difference domain decomposition schemes for
  nonstationary convection-diffusion problems, Differ. Equations 32~(7) (1996)
  929--933.

\bibitem{0971.65076}
A.~A. Samarskii, The theory of difference schemes, Pure and Applied
  Mathematics, Marcel Dekker. 240. New York, NY: Marcel Dekker. 786 p., 2001.

\bibitem{Marchuk:1990:SAD}
G.~I. Marchuk, Splitting and alternating direction methods, in: P.~G. Ciarlet,
  J.-L. Lions (Eds.), Handbook of Numerical Analysis, Vol. I, North-Holland,
  1990, pp. 197--462.

\bibitem{0209.47103}
N.~N. Yanenko, The method of fractional steps. The solution of problems of
  mathematical physics in several variables, Berlin-Heidelberg-New York:
  Springer Verlag, VIII, 160 p. with 15 fig., 1971.

\bibitem{0838.65086}
P.~N. Vabishchevich, Regionally additive difference schemes with a stabilizing
  correction for parabolic problems., Comput. Math. Math. Phys. 34~(12) (1994)
  1573--1581.

\bibitem{0297.35037}
D.~G. Gordeziani, G.~V. Meladze, Simulation of the third boundary value problem
  for multidimensional parabolic equations in an arbitrary domain by
  one-dimensional equations, U.S.S.R. Comput. Math. Math. Phys. 14(1974)~(1)
  (1975) 249--253.

\bibitem{0986.65510}
P.~N. Vabishchevich, V.~A. Verakhovskij, Difference schemes for component-wise
  splitting-decomposition of a domain, Mosc. Univ. Comput. Math. Cybern.
  1994~(3) (1994) 7--11.

\bibitem{dryja2007}
M.~Dryja, X.~Tu, A domain decomposition discretization of parabolic problems,
  Numerische Mathematik 107 (2007) 625--640.

\bibitem{samarskii1998regularized}
A.~A. Samarskii, P.~N. Vabishchevich, Regularized additive full approximation
  schemes, Doklady. Mathematics 57~(1) (1998) 83--86.

\bibitem{0712.65089}
V.~N. Abrashin, A variant of the method of variable directions for the solution
  of multi- dimensional problems of mathematical-physics. i., Differ. Equations
  26~(2) (1990) 243--250.

\bibitem{vabishchevich1996vector}
P.~N. Vabishchevich, Vector additive difference schemes for first-order
  evolutionary equations, Computational mathematics and mathematical physics
  36~(3) (1996) 317--322.

\bibitem{samarskii1992regularized}
A.~A. Samarskii, P.~N. Vabishchevich, Regularized difference schemes for
  evolutionary second order equations, Math. Models and Methods in Applied
  Sciences 2~(3) (1992) 295--315.

\bibitem{abrashin1998numerical}
V.~N. Abrashin, P.~N. Vabishchevich, Vector additive schemes for second-order
  evolution equations, Differential Equations 34~(12) (1998) 1673--1681.

\bibitem{0863.65056}
A.~A. Samarskii, P.~N. Vabishchevich, Vector additive schemes of domain
  decomposition for parabolic problems, Differ. Equations 31~(9) (1995)
  1522--1528.

\bibitem{0965.65119}
P.~N. Vabishchevich, V.~A. Verakhovskij, Domain decomposition vector schemes
  for second-order evolution equations, Mosc. Univ. Comput. Math. Cybern.
  1998~(2) (1998) 1--8.

\bibitem{1156.65084}
P.~N. Vabishchevich, Domain decomposition methods with overlapping subdomains
  for the time-dependent problems of mathematical physics., Comput. Methods
  Appl. Math. 8~(4) (2008) 393--405.

\bibitem{1013.65106}
Y.~Zhuang, X.-H. Sun, Stabilized explicit-implicit domain decomposition methods
  for the numerical solution of parabolic equations., SIAM J. Sci. Comput.
  24~(1) (2002) 335--358.

\bibitem{1120.65098}
Z.~Yu, An alternating explicit-implicit domain decomposition method for the
  parallel solution of parabolic equations., J. Comput. Appl. Math. 206~(1)
  (2007) 549--566.

\bibitem{1076.65079}
T.~Sun, Stability and error analysis on partially implicit schemes., Numer.
  Methods Partial Differ. Equations 21~(4) (2005) 843--858.

\bibitem{1099.65084}
Y.~Jun, T.-Z. Mai, Adi method -- domain decomposition., Appl. Numer. Math.
  56~(8) (2006) 1092--1107.

\bibitem{1094.65097}
Y.~Jun, T.-Z. Mai, Ipic domain decomposition algorithm for parabolic problems.,
  Appl. Math. Comput. 177~(1) (2006) 352--364.

\bibitem{1168.65387}
Y.~Jun, T.-Z. Mai, Numerical analysis of the rectangular domain decomposition
  method., Commun. Numer. Methods Eng. 25~(7) (2009) 810--826.

\bibitem{samarskii1989numerical}
A.~A. Samarskii, E.~S. Nikolaev, Numerical methods for grid equations,
  Birkh{\"a}user, 1989.

\bibitem{pearac}
D.~W. Peaceman, H.~H. Rachford, The numerical solution of parabolic and
  elliptic differential equations, J. SIAM 3 (1955) 28--41.

\bibitem{dourac}
J.~J. Douglas, H.~H. Rachford, On the numerical solution of heat conduction
  problems in two and three space variables, Trans. Amer. Math. Soc. 82 (1956)
  421--439.

\bibitem{samgul}
A.~A. Samarskii, A.~V. Gulin, Stability of difference schemes, Moscow, 1973, in
  Russian.

\bibitem{0184.38503}
G.~Strang, On the construstion and comparison of difference schemes., SIAM J.
  Numer. Anal. 5 (1968) 506--517.

\bibitem{fryazinov1968economical}
I.~V. Fryazinov, The economical symmetrized schemes for the solution of
  boundary value problems for multidimensional parabolic equation, Zh. Vychisl.
  Mat. i Mat. Fiz 8~(2) (1968) 436--443.

\end{thebibliography}

\end{document}